\documentclass[11pt,a4paper]{article}  
\usepackage[left=30mm,right=30mm,top=30mm,bottom=30mm,footskip=12mm]{geometry}

\title{\newheadfont\vskip-1.0em On amenability constants of Fourier algebras: new bounds and new examples}
\author{\sc Y. Choi, M. Ghandehari}
\date{10th February 2026}

\usepackage[dvipsnames]{xcolor}
\usepackage[colorlinks=true,linkcolor=blue,citecolor=red]{hyperref}

\usepackage{amsmath,amssymb,amsfonts}
\numberwithin{equation}{section}

\usepackage{tikz-cd}

\newcommand{\dt}[1]{\textit{#1}} 
\newcommand{\defeq}{\mathbin{:=}}

\newcommand{\bbC}{{\mathbb C}}

\newcommand{\bbS}{{\mathbb S}}
\newcommand{\bbZ}{{\mathbb Z}}

\newcommand{\cU}{{\mathcal U}} 
\newcommand{\cG}{{\mathcal G}}

\newcommand{\cF}{{\mathcal F}}

\newcommand{\FA}{{\rm A}} 

\newcommand{\AM}{{\rm AM}} 
\newcommand{\AD}{{\rm AD}} 

\newcommand{\norm}[1]{{\lVert#1\rVert}}
\newcommand{\ptp}{\widehat{\otimes}}

\newcommand{\Ghat}{\widehat{G}}

\DeclareMathOperator{\maxdeg}{maxdeg}
\DeclareMathOperator{\supp}{supp}
\DeclareMathOperator{\Tr}{Tr}

\newcommand{\FS}{{\rm B}} 
\newcommand{\cL}{{\mathcal L}} 

\newcommand{\sH}{\mathsf{H}} 

\newcommand{\Zadic}[1]{{\mathbb Z}_{#1}} 

\newcommand{\Bdd}{{\mathcal B}} 
\newcommand{\cS}{{\mathcal S}}
\newcommand{\cV}{{\mathcal V}}

\newcommand{\bbN}{{\mathbb N}}

\newcommand{\Hpred}{{\mathsf{Hr}}_p}

\newcommand{\Cst}{{\rm C}^\ast}

\newcommand{\Nhat}{\widehat{N}} 

\newcommand{\utmat}[3]{{\begin{pmatrix} 1 & {#1} & {#2} \\ 0 & 1 & {#3} \\ 0 & 0 & 1 \end{pmatrix}}}

\newenvironment{alphnum}{%
\begin{enumerate}

}{\end{enumerate}\ignorespacesafterend}

\newcommand{\Nat}{{\mathbb N}}
\newcommand{\Real}{{\mathbb R}}
\newcommand{\Cplx}{{\mathbb C}}
\newcommand{\Zahl}{{\mathbb Z}}

\usepackage{amsthm}
\newcounter{pulse}[section]
\numberwithin{pulse}{section}

\newcommand{\newheadfont}{\bfseries} 

\newtheoremstyle{newplain} 
  {\topsep}   
  {\topsep}   
  {\itshape}  
  {0pt}       
  {\newheadfont} 
  {.}         
  {5pt plus 1pt minus 1pt} 
  {}          

\newtheoremstyle{newdef} 
  {\topsep}   
  {\topsep}   
  {\normalfont}  
  {0pt}       
  {\newheadfont} 
  {.}         
  {5pt plus 1pt minus 1pt} 
  {}          

\newtheoremstyle{newrem}
  {\topsep}   
  {\topsep}   
  {\normalfont}  
  {0pt}       
  {\newheadfont} 
  {.}         
  {5pt plus 1pt minus 1pt} 
  {}          

\theoremstyle{newplain}
\newtheorem{thm}[pulse]{Theorem}
\newtheorem{prop}[pulse]{Proposition}
\newtheorem{lem}[pulse]{Lemma}
\newtheorem{cor}[pulse]{Corollary}
\newtheorem{conj}[pulse]{Conjecture}  
\theoremstyle{newdef}
\newtheorem{dfn}[pulse]{Definition}

\theoremstyle{newrem}
\newtheorem{rem}[pulse]{Remark}
\newtheorem{qn}[pulse]{Question} 

\usepackage{sectsty}
\allsectionsfont{\sffamily}
\renewcommand{\newheadfont}{\sffamily\bfseries}

\renewcommand{\dt}[1]{\textcolor{Bittersweet}{\textbf{#1}}} 

\begin{document}

\maketitle

\begin{abstract}
Let $G$ be a locally compact group. If $G$ is finite then the amenability constant of its Fourier algebra, denoted by ${\rm AM}({\rm A}(G))$, admits an explicit formula [Johnson, JLMS 1994];
if $G$ is infinite then no such formula for ${\rm AM}({\rm A}(G))$ is known, although lower and upper bounds were established by Runde [PAMS 2006].
Using non-abelian Fourier analysis, we obtain a sharper upper bound for ${\rm AM}({\rm A}(G))$ when $G$ is discrete. 
Combining this with previous work of the first author [Choi, IMRN 2023], we exhibit new examples of discrete groups and compact groups where ${\rm AM}({\rm A}(G))$ can be calculated explicitly; previously this was only known for groups that are products of finite groups with ``degenerate'' cases.
Our new examples also provide additional evidence to support the conjecture that Runde's lower bound for the amenability constant is in fact an equality.

\bigskip\noindent
MSC2020: 43A30, 46H20 (primary); 20C15, 22D10 (secondary)

\end{abstract}

%

\section{Introduction}

\subsection{Background context}
The Fourier algebra of a locally compact group $G$, denoted by $\FA(G)$, is a Banach algebra of continuous functions on $G$, whose norm encodes the unitary representation theory of~$G$. Unlike the commutative $\Cst$-algebra $C_0(G)$, which only remembers the underlying topological space of~$G$, the Fourier algebra $\FA(G)$ remembers both the underlying topological space of $G$ and the group structure. More precisely, if $G_1$ and $G_2$ are locally compact groups, then the Banach algebras $\FA(G_1)$ and $\FA(G_2)$ are isometrically isomorphic if and only if $G_1$ and $G_2$ are isomorphic as topological groups; this is a theorem of Walter \cite{Walter-1972}. One is therefore led to study how structural properties of a group are reflected in Banach-algebraic properties of its Fourier algebra.

One important structural property of Banach algebras is the notion of \emph{amenability}, whose study has inspired some deep work within particular classes of Banach algebras such as $L^1$-(semi)group algebras or $\Cst$-algebras.
For instance, it was shown in the early 1970s that the $L^1$-convolution algebra of a locally compact group $G$ is amenable if and only if $G$ is amenable; and it was shown in the 1970s-1980s, through the combined work of many authors, that a $\Cst$-algebra is amenable precisely when it is nuclear.
(For self-contained presentations of both results, a good recent source is the book \cite{Runde-newbook}.)
The corresponding characterization for Fourier algebras proved to be more elusive, but was finally obtained by Forrest and Runde in a 2005 paper~\cite{Forrest-Runde-2005}: they proved that the Fourier algebra of a locally compact group $G$ is amenable if and only if $G$ is \dt{virtually abelian} (that is, has an abelian subgroup of finite index).

Within the class of amenable Banach algebras, one can hope to capture more information by considering a quantitative variant. This leads to the notion of the \emph{amenability constant} of a Banach algebra $A$, which we denote by $\AM(A)$. A formal definition will be given in Section~\ref{ss:AM defn}; for the current discussion, we merely note the following properties.
\begin{itemize}
\item $\AM(A)\in [1,\infty]$;
\item $A$ is amenable if and only if $\AM(A)<\infty$;
\item if $A$ and $B$ are isometrically isomorphic, then $\AM(A)=\AM(B)$.
\end{itemize}

While there is a large body of work during the last 50 years that studies amenability for Banach algebras, there has been relatively little investigation of amenability constants. One possible reason is that for $L^1$-group algebras and $\Cst$-algebras, there is a dichotomy: if $A$ is a Banach algebra from one of these two classes, then $\AM(A)$ is either $1$ or~$\infty$, and so the amenability constant is not capturing any further information.
The situation for Fourier algebras is very different, as Johnson showed in \cite[Section~4]{Johnson-1994}, and thirty years after his paper many fundamental questions about the amenability constants of Fourier algebras remain open.
In particular, the following question represents a significant gap in our current state of knowledge.

\begin{qn}\label{q:is AMAG multiplicative}
Let $G_1$ and $G_2$ be locally compact groups. Is $\AM(\FA(G_1\times G_2))$  equal to $\AM(\FA(G_1))\ \AM(\FA(G_2))$?
\end{qn}

To date, no counterexamples are known, but there is no general principle that would imply a positive answer. As we will see in Section~\ref{ss:intro bounds}, the answer is positive if both $G_1$ and $G_2$ are finite, but the proof breaks down as soon as one of the groups is infinite.
The answer is also positive if either $G_1$ or $G_2$ is abelian:
this result was probably known to Johnson, and was certainly known to the authors of \cite{Forrest-Runde-2005}, although we did not find it stated explicitly in the literature. (See Proposition \ref{p:multiply by LCA} for a proof.)

A positive answer to Question~\ref{q:is AMAG multiplicative} would follow from a conjecture of the first author that was implicitly made in \cite{Choi-minorant}. The precise statement will be given as Conjecture \ref{conj:main}; in order to motivate the conjecture, we need to review some results from the papers \cite{Johnson-1994} and~\cite{Runde-2006}.

\subsection{Bounds for the amenability constant of the Fourier algebra}\label{ss:intro bounds}
The amenability constant of the Fourier algebra was first studied by Johnson in~\cite{Johnson-1994}. The present paper is in some sense complementary to Johnson's, since we shall focus on discrete groups while \cite{Johnson-1994} deals exclusively with compact groups.

One of Johnson's discoveries was that when $G$ is \emph{finite}, there is a remarkable formula for the amenability constant of its Fourier algebra:
\begin{equation}\label{eq:johnson formula}
\AM(\FA(G)) = \frac{1}{|G|} \sum_{\pi\in \widehat{G}} (d_\pi)^3
\end{equation}
where $\widehat{G}$ denotes the set of irreducible representations of $G$ (up to equivalence) and $d_\pi$ denotes the degree of the representation~$\pi$. See \cite[Theorem 4.1]{Johnson-1994} for a proof.

It follows easily from \eqref{eq:johnson formula} that Question~\ref{q:is AMAG multiplicative} has a positive answer for finite groups, and as an application of this, Johnson was able to provide explicit examples of compact groups $G$ for which $\FA(G)$ is \emph{non-amenable}.
However, he did not discuss the problem of calculating the exact value of $\AM(\FA(G))$ for infinite non-abelian groups. This seems to be a challenging problem, because we currently do not have any analogue of \eqref{eq:johnson formula} for infinite groups.
Nevertheless, approximately ten years later, Runde showed in \cite{Runde-2006} that one can give useful upper and lower bounds on $\AM(\FA(G))$, in terms of representation-theoretic data associated to~$G$. Namely,
\begin{equation}\label{eq:runde bounds}
1\leq \AD(G) \leq \AM(\FA(G))\leq \maxdeg(G),
\end{equation}
where $\maxdeg(G)$ denotes the supremum of degrees of irreducible unitary representations of~$G$, and $\AD(G)$ denotes\footnotemark\ the Fourier--Stieltjes norm of the anti-diagonal subset of $G\times G$. (Further details will be given in Definition \ref{d:AD(G)}.)
The upper bound is proved in \cite[Lemma 2.7]{Runde-2006}, while the lower bound is proved in \cite[Lemma~3.1]{Runde-2006}.
\footnotetext{This notation is not used in Runde's paper. It is taken from \cite{Choi-minorant}, with the caveat that there $\AD$ was only considered for discrete groups. What we denote by $\AD(G)$ in this paper corresponds to $\AD(G_d)$ in the notation of \cite{Choi-minorant}.}

Runde's proof of the inequality $\AM(\FA(G))\leq\maxdeg(G)$ requires techniques from the theory of operator spaces. However, as he himself observed, for finite groups this inequality is an elementary consequence of Equation~\eqref{eq:johnson formula}. Indeed, for finite $G$, one can actually obtain a \emph{sharper} upper bound
\begin{equation}\label{eq:finite group improved bound}
\AM(\FA(G))
 \leq  1+ (\maxdeg(G)-1) \left( 1- \frac{1}{|[G,G]|}\right)  \,,
 \end{equation}
where $[G,G]$ denotes the commutator subgroup of~$G$.
(The proof just uses basic facts from the character theory of finite groups; details of the calculation will be given in Section~\ref{s:soft results}.)
Furthermore: if $G$ is finite and non-abelian, the right-hand side of \eqref{eq:finite group improved bound} is \emph{strictly} less than $\maxdeg(G)$, so that for such groups the upper bound in \eqref{eq:runde bounds} is \emph{never sharp}. We will return to this point in the next section.

Returning to the setting of general groups: it is not immediately clear what we gain from the lower bound in \eqref{eq:runde bounds}, since $\AD(G)$ seems just as hard to calculate as $\AM(\FA(G))$. However, results of the first author in \cite{Choi-minorant} indicate that one can say quite a lot about $\AD(G)$ for particular examples.
By \cite[Proposition 3.2]{Runde-2006}, $\AD(G)$ is finite if and only if $G$ is virtually abelian, and so we can restrict attention to this class of groups.
It turns out that for countable virtually abelian $G$, there is an explicit formula for $\AD(G)$ in terms of the Plancherel measure on $\widehat{G}$ (see \cite[Theorem 1.5]{Choi-minorant}).

When $G$ is finite, it turns out that $\AD(G)=\AM(\FA(G))$ (\cite[Theorem 1.4]{Choi-minorant}). This was unexpected: the proof in \cite{Runde-2006} of the inequality $\AD(G)\leq\AM(\FA(G))$ gives no reason to suppose that equality should hold for finite groups. For, as discussed in \cite[Section~6]{Choi-minorant}, the inequality arises from considering a norm-decreasing map between two Banach spaces which is known to be non-isometric for every non-abelian group, even the finite ones.

\subsection{New results}\label{ss:intro new results}

It was asked in~\cite[Question 6.1]{Choi-minorant} if the inequality $\AD(G)\leq\AM(\FA(G))$ is actually an equality in all cases. We now state this as a formal conjecture.

\begin{conj}\label{conj:main}
$\AM(\FA(G))=\AD(G)$ for every locally compact group $G$.
\end{conj}

By \cite[Proposition 2.9]{Choi-minorant}, $\AD$ is ``multiplicative'' with respect to direct products of groups. Hence Conjecture \ref{conj:main} implies that Question~\ref{q:is AMAG multiplicative} has a positive answer.
It also follows from \cite[Proposition 2.9]{Choi-minorant} that if $\theta:H\to G$ is any continuous injective group homomorphism, then $\AD(H)\leq \AD(G)$, and so Conjecture \ref{conj:main} implies that $\AM(\FA(H))\leq\AM(\FA(G))$. These improved functorial properties of $\AM(\FA(\cdot))$ would then provide powerful additional tools to aid in calculating or estimating amenability constants of Fourier algebras.

Prior to the present paper, the reasons to believe Conjecture \ref{conj:main} were more theoretical than empirical. Indeed, the conjecture was only known to hold in the following cases.
\begin{enumerate}
\item
If $G$ is not virtually abelian, then it satisfies Conjecture \ref{conj:main}, but for a somewhat unsatisfying reason: namely, for such $G$, we have $\AD(G)=+\infty$. This follows from the results in \cite{Forrest-Runde-2005}.
\item
If $G$ is finite, then it satisfies Conjecture \ref{conj:main} (as already discussed).
\item\label{li:AD-maximal}
If $G$ is virtually abelian and $\AD(G)=\maxdeg(G)$, then it satisfies Conjecture \ref{conj:main}. (This is immediate from the inequalities in~\eqref{eq:runde bounds}.) For convenience, we shall call groups with these properties \dt{AD-maximal}.
\end{enumerate}

The main purpose of this paper is to obtain new examples of groups $G$ for which $\AM(\FA(G))$ can be calculated explicitly, which then provide further evidence for Conjecture~\ref{conj:main}. 
We start by presenting some new results that can be proved with soft techniques.

\begin{thm}[Countable saturation]
\label{t:AMA-discrete-countable}
Let $G$ be a discrete group. Then there is a countable subgroup $\Gamma\leq G$ such that $\AM(\FA(G))=\AM(\FA(\Gamma))$.
\end{thm}

Theorem \ref{t:AMA-discrete-countable} is slightly counterintuitive, because when $\Gamma$ is a subgroup of $G$ every bounded approximate diagonal for $\FA(\Gamma)$ annihilates any element of $\FA(G)$ whose support is disjoint from~$\Gamma$.

Combining Theorem \ref{t:AMA-discrete-countable} with known hereditary properties of $\AD$ yields the following result.

\begin{cor}\label{c:conjecture-discrete-countable}
If Conjecture \ref{conj:main} holds for all countable groups, then it holds for all discrete groups.
\end{cor}

As we already mentioned, $\AD$ is ``multiplicative'' with respect to direct products of groups. It follows that the class of groups satisfying Conjecture \ref{conj:main} is closed under taking finite products (Proposition \ref{p:taking products preserves conjecture}). Applying this to the known examples that were listed earlier yields the following result.

\begin{cor}\label{c:cheap cases}
Let $F$ be a finite group and let $H$ be an AD-maximal group. Then $F\times H$ satisfies Conjecture~\ref{conj:main}.
\end{cor}

\begin{rem}\label{r:products of AD-maximal}
It is tempting to state a more general version of Corollary \ref{c:cheap cases}, where one takes products of finitely many groups, each of which are either finite or AD-maximal. However, this would not lead to a larger class of groups, because of the following fact: \emph{the product of two AD-maximal groups is AD-maximal}. This fact is a straightforward consequence of known results; details are given in Appendix~\ref{app:maxdeg}.
\end{rem}

It is natural to ask how wide the class of groups covered by Corollary \ref{c:cheap cases}~is. It is clear from the definition that all locally compact abelian groups are AD-maximal, but not clear that there are any other examples. However, \cite[Theorem 1.10]{Choi-minorant} shows how to construct non-abelian examples of AD-maximal groups that are either discrete, non-connected Lie, or compact.

While Corollary \ref{c:cheap cases} provides new examples to confirm Conjecture \ref{conj:main}, these examples may seem somewhat artificial. 
The next result, which is one of the main results of this paper, shows that we can find natural examples of groups that satisfy  Conjecture \ref{conj:main} and are not covered by Corollary \ref{c:cheap cases}.

\begin{thm}\label{t:new examples}
Let $p\geq 2$ be prime.
There exist a discrete group $\Hpred(\bbZ)$ and a compact group $\Hpred(\Zadic{p})$ with the following properties.
\begin{alphnum}
\item\label{li:p-reduced integer Heis}
$\AM(\FA(\Hpred(\bbZ)))=\AD(\Hpred(\bbZ))=p-1+p^{-1}$.
\item\label{li:p-reduced p-adic Heis}
$\AM(\FA(\Hpred(\Zadic{p})))=\AD(\Hpred(\Zadic{p}))=p-1+p^{-1}$.
\item\label{li:not cheap}
Neither $\Hpred(\bbZ)$ nor $\Hpred(\Zadic{p})$ are isomorphic to groups of the form described in Corollary~\ref{c:cheap cases}.
\end{alphnum}
\end{thm}

The groups $\Hpred(\bbZ)$ and $\Hpred(\Zadic{p})$ arise as quotients of Heisenberg groups over the integers and $p$-adic integers, respectively; details are given in Section~\ref{s:p-reduced Heis}.

In fact, part \ref{li:p-reduced integer Heis} of Theorem \ref{t:new examples} is a special case of the following new result.

\begin{thm}[Groups with two degrees of irreducible representations]\label{t:two degrees}
Let $d\in\Nat$. Suppose that $G$ is a countable group such that $\{ d_\pi \colon \pi\in \widehat{G}\} = \{1,d\}$.
Then
\[
\AM(\FA(G)) = \AD(G) = 1 + (d-1)\left( 1- \frac{1}{|[G,G]|} \right) \,.
\]
\end{thm}

The proof of Theorem \ref{t:two degrees} has two separate ingredients. The first ingredient is the formula for $\AD(G)$ in terms of the Plancherel measure of $G$, obtained in \cite[Theorem~1.5]{Choi-minorant}. The second ingredient, which will be stated later as Theorem~\ref{t:sharper upper bound}, is the fact that when $G$ is discrete and virtually abelian then $\AM(\FA(G))$ satisfies the same sharper upper bound that was described in \eqref{eq:finite group improved bound}.
This is a surprising result, since the existing proof of \eqref{eq:finite group improved bound} is based on Johnson's formula \eqref{eq:johnson formula} which is not available for infinite groups. Instead, our proof uses the noncommutative Fourier transform and inverse Fourier transform for countable virtually abelian groups, together with a bootstrapping argument that reduces the general case to the countable case.

Our final new result is another application of Theorem \ref{t:sharper upper bound}, motivated by \cite[Question~6.2]{Choi-minorant}.
To make our narrative more self-contained, we combine the relevant results from the older paper with the new result, and state the combination as the following theorem.

\begin{thm}
\label{t:minimal AM(A(non-abelian))}
Let $G$ be a locally compact \emph{non-abelian} group.
\begin{alphnum}
\item\label{li:AD geq 3/2}
$\AM(\FA(G))\geq\AD(G) \geq 3/2$.
\item\label{li:minimal AD}
$\AD(G)=3/2$ if and only if $|G\mathbin{:}Z(G)|=4$.
\item\label{li:partial converse}
If $G$ is discrete and $|G\mathbin{:}Z(G)|=4$, then $\AM(\FA(G))=3/2$.
\end{alphnum}
\end{thm}

Part \ref{li:AD geq 3/2} is \cite[Theorem 1.6]{Choi-minorant} and part \ref{li:minimal AD} is \cite[Theorem 1.7]{Choi-minorant}; it is part \ref{li:partial converse} that is new.
Taken together, the three parts of this theorem characterize those discrete non-abelian groups for which $\AM(\FA(G))$ attains its minimal value, and provide a positive answer to \cite[Question 6.2]{Choi-minorant} in the special case of discrete groups.

We now summarize the structure of this paper.
Since it is aimed at a general audience, we shall start by reviewing the necessary background material on amenability constants, Fourier algebras, and the operator-valued Fourier transform; this takes up most of Section~\ref{s:prelim}.

In Section \ref{s:soft results} we give proofs for the ``soft'' results that were stated earlier in this section; these are proved by general functional-analytic principles.
Section \ref{s:prove upper bound} is the heart of the paper, in which we state and prove the sharper upper bound on $\AM(\FA(G))$ when $G$ is discrete, and use it to prove Theorems \ref{t:two degrees} and \ref{t:minimal AM(A(non-abelian))}. In Section \ref{s:p-reduced Heis}, we draw on results proved in the earlier sections to show that the groups $\Hpred(\bbZ)$ and $\Hpred(\Zadic{p})$ satisfy the properties listed in Theorem~\ref{t:new examples}. We finish with some concluding remarks and suggestions for future work, in Section~\ref{s:where next}.

Some results on which we rely appear to be folklore: they are known to experts, but we were unable to find references in the literature that stated these results in the form that we needed. To make the paper more self-contained, we have included some of the relevant details and explanations of these folklore results in an appendix.

\section{Preliminaries}\label{s:prelim}

\subsection{Notation and other conventions}
We start by fixing some conventions and notation.
All vector spaces and algebras are over the complex field. The tensor product of two vector spaces $V$ and $W$ is denoted by $V\otimes W$; if $V$ and $W$ are Banach spaces, then $V\ptp W$ denotes their projective tensor product. If $A$ is a Banach algebra, then the multiplication map on $A$ extends to a contractive linear map $\Delta: A\ptp A \to A$.

Given a Hausdorff topological space $X$ and a continuous function $f:X\to\bbC$, we denote the support of $f$ by $\supp(f)$. 
If $X$ is a locally compact Hausdorff space, we write $C_c(X)$ for the set of compactly supported functions $X\to\Cplx$. If $X$ is discrete, this space is the same as the space of finitely supported functions $X\to \Cplx$, and in this setting we shall use the notation $c_{00}(X)$.

By a \dt{directed set}, we mean a non-empty partially ordered set $(I,\preceq)$ with the following additional property: for each $i,j\in I$, there exists $k\in I$ such that $i\preceq k$ and $j\preceq k$. A subset $I'\subseteq I$ is said to be \dt{cofinal} (with respect to $\preceq$) if, for each $j\in I$, there exists some $k\in I'$ such that $j\preceq k$.

Given a directed set $(I,\preceq)$ and some set~$S$, a \dt{directed system of subsets of $S$} is a family of sets  $(S_i)_{i\in I}$ such that $S_i\subseteq S_j\subseteq S$ for all $i\preceq j$. One defines a directed system of subgroups of some given group~$G$, or a directed system of (closed) subalgebras of some given (Banach) algebra~$A$, in the same way.

To reduce repetition: all representations of topological groups on Hilbert spaces are assumed to be unitary and SOT-continuous.
As a further notational simplification:
given such a representation $\pi:G\to \cU(H_\pi)$ and $p\in [1,\infty)$, we write $S^p(\pi)$ for the space of Schatten-$p$ operators on the Hilbert space~$H_\pi$.

\subsection{Amenability constants of Banach algebras}
\label{ss:AM defn}

The following is not Johnson's original definition of amenability for Banach algebras, but was shown by him in \cite{Johnson-1972_appdiag} to be an equivalent definition.

\begin{dfn}\label{d:AM(A)}
Let $A$ be a Banach algebra. A \dt{bounded approximate diagonal} for $A$ (b.a.d.~for short) is a bounded net $(m_\alpha)$ in $A\ptp A$ with the following two properties:
\begin{alphnum}
\item\label{li:app-central}
for each $a\in A$, $\norm{a\cdot m_\alpha- m_\alpha\cdot a}\to 0$;
\item\label{li:app-id}
for each $a\in A$, $\norm{a\Delta(m_\alpha)-a}\to 0$.
\end{alphnum}
If $A$ has a b.a.d.\ we say that $A$ is \dt{amenable}. The \dt{amenability constant} of $A$, denoted by $\AM(A)$, is defined to be the infimum of $\sup_\alpha\norm{m_\alpha}$ over all choices of b.a.d.\ for~$A$, with the convention that $\AM(A)=+\infty$ when $A$ is non-amenable.
\end{dfn}

\begin{rem}\label{r:inf attained}
Taking a weak$^*$-cluster point of a b.a.d.\ for $A$, inside $(A\ptp A)^{**}$, leads to the concept of a \dt{virtual diagonal} for~$A$, and one can equivalently define $\AM(A)$ to be the infimum of norms of all possible virtual diagonals for~$A$. In fact, since the set of virtual diagonals for $A$ is a weak$^*$-closed subset of $(A\ptp A)^{**}$, the proof of \cite[Lemma 1.2]{Johnson-1972_appdiag} shows that whenever $A$ is amenable, there exists a b.a.d.\ for $A$ whose norm is \emph{equal} to $\AM(A)$.
\end{rem}

\subsection{General properties of Fourier algebras}
As in the introduction, the Fourier algebra of a locally compact group $G$ is denoted by $\FA(G)$.
Most of what we need can be found in the original paper of Eymard \cite{Eymard-1964}, but a convenient reference that we shall use for specific details is the book \cite{Kaniuth-Lau-book}.

There are several different ways to define $\FA(G)$. All of them give the same object in the end, but in any given situation, one definition may be more convenient than the others.
For this paper, the most natural definition (not the original one) is as follows: fix a Haar measure on $G$, and let $\lambda:G\to \cU(L^2(G))$ be the left regular representation. There is a bounded bilinear map $L^2(G)\times L^2(G) \to C_0(G)$, which sends a pair $(\xi,\eta)$ to the \dt{coefficient function}
\[
s\mapsto \langle \lambda(s)\xi, \overline{\eta} \rangle = \int_G \xi(s^{-1}t)\eta(t)\,dt\,.
\]
Therefore, by the universal property of the projective tensor product, we obtain a contractive linear map $\alpha_G: L^2(G)\ptp L^2(G) \to C_0(G)$ that satisfies $\alpha_G(\xi\otimes\overline{\eta}) = \langle\lambda(\cdot)\xi, \eta\rangle$. We now define $\FA(G)$ to be the image of $\alpha_G$, equipped with the \emph{quotient norm} inherited from $(L^2(G)\ptp L^2(G))/\ker\alpha_G$.

\begin{rem}\label{r:standard form}
It follows from deep results about group von Neumann algebras that each $f\in\FA(G)$ can be realized as a single coefficient function $f=\alpha_G(\xi\otimes\eta)$ for some $\xi,\eta\in L^2(G)$ that satisfy $\norm{f}_{\FA(G)} = \norm{\xi}\,\norm{\eta}$. (See e.g.~\cite[Theorem 2.4.3]{Kaniuth-Lau-book}.) This is sometimes convenient for calculations, since it bypasses the need to consider linear combinations of elementary tensors.
\end{rem}

Another result about Fourier algebras is sufficiently useful that we state it here for emphasis, and also to set up some notation.

\begin{lem}\label{l:embed and restrict}
If $G$ is a discrete group and $H$ is a subgroup, and $f:H\to\Cplx$, define $\iota_H f : G\to \Cplx$~by
\[
\iota_H f (x) = \begin{cases} f(x) & \text{if $x\in H$} \\ 0 & \text{otherwise.} \end{cases}
\]
Then the following properties hold:
\begin{alphnum}
\item\label{li:embedding}
$\iota_H$ is an isometric algebra homomorphism from $\FA(H)$ into $\FA(G)$;
\item\label{li:restriction}
restriction of functions from $G$ to $H$ defines a contractive algebra homomorphism $\rho_H: \FA(G)\to\FA(H)$;
\item\label{li:range of embedding}
$\iota_H(\FA(H)) = \{ f\in \FA(G) \colon \supp(f)\subseteq H\}$.
\end{alphnum}
\end{lem}

\begin{proof}
Parts \ref{li:embedding} and \ref{li:restriction} are special cases of more general results for open subgroups of locally compact groups: see e.g. \cite[Pro\-pos\-ition 2.4.1]{Kaniuth-Lau-book}.
Part \ref{li:range of embedding} is a straightforward consequence of \ref{li:embedding} and \ref{li:restriction}; details are left to the reader.
\end{proof}

\begin{rem}\label{r:Herz restriction}
Lemma \ref{l:embed and restrict}\ref{li:restriction} can be generalized significantly, although the proof is substantially harder. Namely: for any locally compact group $G$ and any closed subgroup $H$, restriction of functions from $G$ to $H$ defines a contractive and surjective algebra homomorphism $\rho_H:\FA(G)\to\FA(H)$. This is commonly known as \dt{Herz's restriction theorem}: see
\cite[Section 2.6]{Kaniuth-Lau-book} for a readable account.
\end{rem}

\subsection{The Fourier anti-diagonal constant of a group}

The lower bound on $\AM(\FA(G))$ in \eqref{eq:runde bounds} is given by an invariant $\AD(G)$; although it is not the main object of interest in this paper, we briefly give its definition and state some of the properties that we will need. For the background context and definitions concerning Fourier--Stieltjes algebras, the reader is referred to \cite[Section 2.1]{Kaniuth-Lau-book}.

Given a locally compact group $G$ we write $G_d$ for the same group but equipped with the discrete topology. We write $\operatorname{adiag}(G)$ for the set $\{(x,x^{-1})\colon x\in G\} \subseteq G\times G$. The key part of Forrest and Runde's proof that amenability of $\FA(G)$ forces $G$ to be virtually abelian, is the (deep) fact that $1_{\operatorname{adiag}(G)}$ belongs to the Fourier--Stieltjes algebra $\FS(G_d\times G_d)$ if and only if $G_d$ is virtually abelian.

\begin{dfn}\label{d:AD(G)}
We define $\AD(G) \defeq \norm{1_{\operatorname{adiag}(G)}}_{\FS(G_d\times G_d)}$, with the convention that this equals infinity if the function does not belong to $\FS(G_d\times G_d)$.
\end{dfn}

Note, in particular, that $\AD(G)$ does not depend on the original topology of~$G$.

\begin{prop}[Functorial properties of $\AD$]\label{p:AD properties}\
\begin{alphnum}
\item\label{li:AD of subgroup} (Monotonicity) Let $G$ be a group and let $H$ be a subgroup of~$G$. Then $\AD(H)\leq\AD(G)$.
\item\label{li:AD of product} (Respects products) Let $G_1$ and $G_2$ be groups. Then $\AD(G_1\times G_2)=\AD(G_1)\AD(G_2)$.
\end{alphnum}
\end{prop}

\begin{proof}
See \cite[Proposition 2.9]{Choi-minorant}.
\end{proof}

\begin{rem}
It is also shown in \cite[Proposition 2.10]{Choi-minorant} that $\AD$ has a ``countable saturation'' property: for any $G$ we can always find a countable subgroup $\Gamma\leq G$ such that $\AD(\Gamma) = \AD(G)$. Although we do not need this result, we note that if Conjecture \ref{conj:main} is true then it would imply that $\AM(\FA(\cdot))$ has the same countable saturation property. Theorem \ref{t:AMA-discrete-countable} is therefore a partial confirmation of this prediction.
\end{rem}

\subsection{The operator-valued Fourier transform}
Given a locally compact group $G$, a 
representation $\pi:G \to \cU(H_\pi)$, and a finite Borel measure $\mu$ on $G$, we define
\[
\pi(\mu) \defeq \int_G \pi(x) d\mu(x) \in \Bdd(H_\pi)
\]
where the integral is understood in the weak sense. In the specific case where $\pi$ is irreducible, we refer to $\pi(\mu)$ as an \dt{operator-valued Fourier coefficient of $\mu$}.
Note that when $G$ is discrete, a finite measure on $G$ is simply given by some $f\in\ell^1(G)$, and we have
$\pi(f) = \sum_{s\in G} f(s)\pi(s)$.

\begin{prop}\label{p:pi of mu_K}
Let $G$ be a locally compact group, let $K$ be a compact normal subgroup of~$G$, and let $\mu_K$ be the pushforward of Haar measure on $K$ along the inclusion map $K \hookrightarrow G$, normalized such that $\mu_K(K)=1$.

Let $\pi: G\to\cU(H_\pi)$ be an irreducible unitary representation, and let $I_\pi$ denote the identity operator on $H_\pi$. If $\pi$ factors through the quotient homomorphism $G\to G/K$, then $\pi(\mu_K)=I_\pi$; if not, then $\pi(\mu_K)=0$.
\end{prop}

\begin{proof}
Suppose that $\pi$ factors through the quotient homomorphism $q: G\to G/K$, i.e., $\pi(k)=I_{\pi}$ for every $k\in K$. Then, for arbitrary $\xi,\eta\in H_\pi$, we have
\[
\langle \pi(\mu_K)\xi,\eta\rangle=\int_G \langle \pi(x)\xi,\eta\rangle\, d\mu_K(x)=\int_K \langle \pi(x)\xi,\eta\rangle\, d\mu_K(x)=\langle \xi,\eta\rangle \int_K \, d\mu_K(x)=\langle \xi,\eta\rangle,
\]
and so $\pi(\mu_K)=I_\pi$.

Now suppose that $\pi$ does not factor through~$q$.
For each $g\in G$, conjugation by $g$ is a continuous automorphism of $K$, so by uniqueness of normalized Haar measure on $K$ we see that $\delta_g * \mu_K * \delta_{g^{-1}}=\mu_K$. Hence, $\pi(g)\pi(\mu_K)\pi(g)^{-1}=\pi(\mu_K)$ for all $g\in G$, which implies by Schur's lemma that $\pi(\mu_K)=\lambda_\pi I_\pi$ for some $\lambda_\pi\in\Cplx$.
Since $\pi$ does not factor through $q$, there exists some $k\in K$ with $\pi(k)\neq I_\pi$. But since $\delta_k*\mu_K=\mu_K$, we have $\lambda_\pi \pi(k)=\lambda_\pi I_\pi$. This forces $\lambda_\pi=0$.
\end{proof}

We now specialize to cases where $G$ is second-countable, unimodular and Type~I. 
Given such a group $G$ and a fixed choice of Haar measure on $G$, 
there exists a unique Radon measure $\nu$ on $\widehat{G}$ (the \dt{Plancherel measure} for $G$) that satisfies
\begin{equation}\label{eq:plancherel formula}
\norm{f}_2^2 = \int_{\widehat{G}} \Tr \left( \pi(f)\pi(f)^* \right) \;d\nu(\pi) \quad(f\in C_c(G)).
\end{equation}
Equation \eqref{eq:plancherel formula} is often referred to as the \dt{Plancherel formula} for $G$; 
for further context and references, see \cite[Section 8.5]{Deitmar-Echterhoff-book}, \cite[Theorem 7.36]{Folland-book} or \cite[Theorem 3.31]{Fuhr-book}.

\begin{rem}\label{r:LCVA implies Type I}
We will apply this machinery in the setting of countable virtually abelian groups. To justify this, we note that every virtually abelian group $G$ has a \emph{normal} abelian subgroup $M$ of finite index (this is a standard observation in group theory), and hence if we regard $G$ as a discrete group, it is Type~I and satisfies $\maxdeg(G)\leq \lvert G\mathbin{:} M\rvert$ (see for example \cite[Satz 5]{Thoma-1968}). 
More generally, if $G$ is locally compact and has an abelian subgroup $H_0$ of finite index, then the closure of $H_0$ in $G$ is an \emph{open} abelian subgroup of finite index. 
This implies that $\maxdeg(G)<\infty$ (see \cite[Propos\-ition~2.1]{Moore-fd-irrep} for a proof using a version of Frobenius reciprocity), which in particular implies that $G$ is Type~I. 
\end{rem}

Later in the paper, to prove Theorem \ref{t:sharper upper bound},
we rely on the fact that when $G$ and $N$ are countable virtually abelian groups one can express the norms on $\FA(G)$, $\FA(G)\ptp\FA(N)$ and $\FA(G\times N)$ in terms of the Fourier transforms for $G$ and~$N$.
Since it is not immediate how to extract these known results from the literature, we shall give a precise statement in Proposition \ref{p:norms via FT} below, and then provide further details/explanation in the appendix (Section~\ref{app:norms via FT}).

Our approach goes via the Plancherel theorem for second-countable unimodular Type~I groups, and so we shall state the required results for this more general class of groups, in case this is useful for future work.
Let $(L^1\cap L^2)(G)$ denote the intersection of $L^1(G)$ and $L^2(G)$, and let
\begin{equation}\label{eq:the nice subspace}
\cV_2^G \defeq \operatorname{span} \{ \alpha_G(\xi\otimes\eta) \colon \xi,\eta\in (L^1\cap L^2)(G)\}
\end{equation}
(this is the uncompleted linear span). Note that $\cV_2^G$ is a dense subspace of $\FA(G)$.

\begin{prop}\label{p:norms via FT}
Let $G$ and $N$ be second-countable unimodular Type~I groups. Fix Haar measures on $G$ and~$N$, and let $\nu_G$ and $\nu_N$ be the corresponding Plancherel measures on $\widehat{G}$ and $\widehat{N}$ respectively.
\begin{alphnum}
\item\label{li:AG-norm}
For every $f\in \cV_2^G$, we have
\begin{equation}\label{eq:A(G)-norm via plancherel}
\norm{f}_{\FA(G)} = \int_{\pi\in \widehat{G}} \norm{\pi(f)}_{S^1(\pi)} \; d\nu_G(\pi)\,.
\end{equation}

\item\label{li:AG-tp-norm}
For every $w\in \cV_2^G \otimes \cV_2^N$, we have
\begin{equation}
\norm{w}_{\FA(G)\ptp\FA(N)}
 = \int_{\widehat{G}\times\widehat{N}} \norm{ (\pi\otimes \sigma)(w)}_{S^1(\pi)\ptp S^1(\sigma)} \;d(\nu_G\times\nu_N)(\pi,\sigma)\,.
\end{equation}

\item\label{li:AGN-norm}
For every $w\in \cV_2^G\otimes \cV_2^N$, we have
\begin{equation}
\norm{w}_{\FA(G\times N)}
 = \int_{\widehat{G}\times\widehat{N}} \norm{ (\pi\otimes \sigma)(w)}_{S^1(\pi\otimes\sigma)} \;d(\nu_G\times\nu_N)(\pi,\sigma) \,.
\end{equation}
\end{alphnum}
\end{prop}

\begin{rem}\label{r:scope}
One can show that part~\ref{li:AG-norm} of Proposition \ref{p:norms via FT} remains valid for all $f\in \FA(G)\cap L^1(G)$,
but doing so requires careful consideration of the pointwise behaviour of the non-abelian Fourier transform and Fourier inversion formula. Similar comments apply to parts \ref{li:AG-tp-norm} and \ref{li:AGN-norm}.
For technical reasons, the space $\cV_2^G$ is easier to handle than $\FA(G)\cap L^1(G)$, and it suffices for our needs.

In fact, we will only apply Proposition \ref{p:norms via FT} when $G$ and  $N$ are countable (in particular, discrete). For such groups, note that each point mass $\delta_s$ arises as $\alpha_G(\delta_e\otimes\delta_s)$, and so by taking linear combinations we have $c_{00}(G)\subseteq \cV_2^G$. Then, upon tensoring, we have
\[
c_{00}(G\times N)= c_{00}(G)\otimes c_{00}(N) \subseteq \cV_2^G\otimes \cV_2^N \;.
\]
Hence, in the intended applications, we can apply Proposition \ref{p:norms via FT} to any $f\in c_{00}(G)$ and any $w\in c_{00}(G\times N)$.
\end{rem}

\section{General results for the amenability constant}
\label{s:soft results}
We start by honouring a promise made in the introduction: namely, we show how to improve on Runde's upper bound $\AM(\FA(G))\leq \maxdeg(G)$ when $G$ is finite.

Let $G$ be a finite group, and write $\cL$ for the set of $1$-dimensional representations of $G$. Using Johnson's formula \eqref{eq:johnson formula}, we have
\[
\begin{aligned}
\AM(\FA(G)) - 1  & = \frac{1}{|G|} \sum_{\pi\in \widehat{G}} (d_\pi)^3 -\frac{1}{|G|} \sum_{\pi\in \widehat{G}} (d_\pi)^2 \\
& = \frac{1}{|G|} \sum_{\pi\in \widehat{G} \setminus\cL} (d_\pi)^3-(d_\pi)^2 \\
& \leq \frac{\maxdeg(G)-1}{|G|} \sum_{\pi\in \widehat{G} \setminus\cL} (d_\pi)^2
& = \frac{\maxdeg(G)-1}{|G|} \left( |G|- |\cL| \right) \,.
\end{aligned}
\]
Let $[G,G]$ denote the commutator subgroup of $G$. Then $G/[G,G]$ is abelian and its Pontrjagin dual can be naturally identified with $\cL$.
In particular, $|\cL| = |G|\,|[G,G]|^{-1}$, and we thus obtain 
\[
\AM(\FA(G))
 \leq  1+ (\maxdeg(G)-1) \left( 1- \frac{1}{|[G,G]|}\right)  \,,
\]
which is precisely the inequality stated in \eqref{eq:finite group improved bound}. As we will see in Section \ref{s:prove upper bound}, the same inequality remains true for infinite discrete groups $G$, but instead of Johnson's formula we have to use norm estimates in tensor products of Fourier algebras.

Next, we give proofs for some properties of $\AM(\FA(\cdot))$ that were stated in the introduction. It is convenient to isolate a common part of these proofs as a separate lemma.

\begin{lem}\label{l:AMAG is submult}
Let $G$ and $H$ be locally compact groups. Then
\[
\AM(\FA(G\times H))\leq \AM(\FA(G))\,\AM(\FA(H)).
\]
\end{lem}

\begin{proof}
There is a natural contractive homomorphism $\FA(G)\ptp\FA(H)\to \FA(G\times H)$, defined by sending an elementary tensor $f\otimes g$ to the function $(s,t)\mapsto f(s)g(t)$.
By standard hereditary properties of amenability constants
(see e.g. \cite[Propositions 2.3.1 and 2.3.14]{Runde-newbook}), it follows that
\[ \AM(\FA(G\times H)) \leq \AM(\FA(G))\, \AM(\FA(H)),
\]
as claimed.
\end{proof}

\begin{prop}\label{p:multiply by LCA}
Let $G$ and $H$ be locally compact groups with $H$ abelian.
Then $\AM(\FA(G))=\AM(\FA(G\times H))$.
\end{prop}

\begin{proof}
There is a contractive algebra homomorphism $\rho_G: \FA(G\times H) \to \FA(G)$, obtained by identifying $G$ with the subgroup $G\times\{e\}\subseteq G\times H$ and appealing to Herz's restriction theorem (see Remark \ref{r:Herz restriction}). Therefore, by a standard property of amenability constants (see e.g.~\cite[Pro\-pos\-ition~2.3.1]{Runde-newbook}) we have $\AM(\FA(G))\leq \AM(\FA(G\times H))$.

We now prove the converse inequality. Since $H$ is abelian, we have $\AM(\FA(H))=1$. (One can deduce this from results concerning operator space tensor products and operator amenability. Alternatively, note that if $\Gamma$ denotes the Pontrjagin dual of $H$, then $\FA(H)$ is isometrically isomorphic to the convolution algebra $L^1(\Gamma)$, and $\AM(L^1(\Gamma))=1$ by \cite[Corollary 1.10]{Stokke-2004}.)
Therefore, Lemma \ref{l:AMAG is submult} gives $\AM(\FA(G\times H))\leq\AM(\FA(G))$, as required.
\end{proof}

\begin{prop}\label{p:taking products preserves conjecture}
Let $G_1,\dots, G_n$ be locally compact virtually abelian groups and let $G=\prod_{i=1}G_i$ (with the product topology). If $\AM(\FA(G_i))=\AD(G_i)$ for all $i$, then $\AM(\FA(G))=\AD(G)$.
\end{prop}

\begin{proof}
We give the proof for $n=2$; the general case follows by induction. 

From Runde's lower bound~\eqref{eq:runde bounds} we know that $\AM(\FA(G_1\times G_2))\geq \AD(G_1\times G_2)$, so it suffices to prove the converse inequality.
Since $G_1$ and $G_2$ satisfy Con\-jecture \ref{conj:main}, and since $\AD$ respects direct products (Proposition \ref{p:AD properties}\ref{li:AD of product}), we have
\[
\AM(\FA(G_1))\,\AM(\FA(G_2)) = \AD(G_1)\, \AD(G_2) = \AD(G_1\times G_2).
\]
Applying Lemma~\ref{l:AMAG is submult} to the left-hand side yields $\AM(\FA(G_1\times G_2)) \leq \AD(G_1\times G_2)$, as required.
\end{proof}

We now prove that when $G$ is discrete, there is always some countable subgroup $\Gamma\leq G$ such that $\AM(\FA(\Gamma))=\AM(\FA(G))$ (Theorem \ref{t:AMA-discrete-countable}).
The key idea is that since $G$ is discrete, we can write $\FA(G)$ as an inductive limit of Fourier algebras of countable subgroups of~$G$. Before spelling this out precisely, we review some general facts concerning amenability constants of inductive limits of Banach algebras.

\begin{lem}[{\cite[Proposition 2.3.15]{Runde-newbook}}]
\label{l:AM of inductive limit}
Let $A$ be a Banach algebra, let $(I,\preceq)$ be a directed set, and let $(A_i)_{i\in I}$ be a directed system of closed subalgebras of $A$.
If 
$A=\overline{\bigcup_{i\in I} A_i}$, then $\AM(A)\leq\sup_{i\in I}\AM(A_i)$. 
\end{lem}

Note that in this result we do not presuppose that $A$ is amenable. Indeed, in many applications, amenability of $A$ is shown to follow from the existence of a suitable ``exhaustion'' of $A$ by a family $(A_i)$ of amenable subalgebras with control on the amenability constant.

\begin{rem}\label{r:limsup}
Given a directed set $(I,\preceq)$ and a function $h: I\to \Real$, we write $\limsup_{i\in I} h(i)$ for $\inf_{i\in I} \sup_{i\preceq j} h(j)$. Since the union $\bigcup_{i\in I} A_i$ is unchanged if we replace $I$ by some cofinal subset $I'$, it follows that $\AM(A)\leq\limsup_{i\in I} \AM(A_i)$.
\end{rem}

The previous lemma and remark are well known to specialists in the amenability of Banach algebras. However, the following refinement seems to be new.

\begin{prop}[Sequential control of the amenability constant of an inductive limit]
\label{p:AM upper bound sequential}
Let $A$ be a Banach algebra. Let $(I,\preceq)$ be a directed set, and let $(A_i)_{i\in I}$ be a directed system of closed subalgebras of $A$, such that $A=\overline{\bigcup_{i\in I} A_i}$. Then there exists an increasing sequence $i(1)\preceq i(2) \preceq \dots$ such that $\AM(A) \leq\limsup_{n\geq 1} \AM(A_{i(n)})$.
\end{prop}

\begin{proof}
We first deal with the case where $\AM(A)<\infty$. By Lemma \ref{l:AM of inductive limit} we can choose $i(1)\in I$ such that $\AM(A_{i(1)}) \geq \AM(A) -2^{-1}$. Now suppose that $n\geq 1$ and we have chosen $i(n)$ such that $\AM(A_{i(n)}) \geq \AM(A)-2^{-n}$. Let $I_n = \{ j\in I \colon j\succeq i(n)\}$. Since $I_n$ is cofinal in $I$, we have $\bigcup_{j\in I_n}A_j = \bigcup_{i\in I} A_i$, so by using Lemma \ref{l:AM of inductive limit} again
we can choose $i(n+1)\succeq i(n)$ such that $\AM(A_{i(n+1)}) \geq \AM(A) - 2^{-(n+1)}$. Continuing by induction, we obtain an increasing sequence $i(1)\preceq i(2)\preceq\dots$ with the desired property.

If $\AM(A)=\infty$ (i.e. if $A$ is not amenable) then the argument is similar, except that our inductive hypothesis has the form $\AM(A_{i(n)}) \geq n$. (The key point is that Lemma \ref{l:AM of inductive limit} remains true if $\AM(A)=\infty$.)
\end{proof}

\begin{rem}
In the previous result, two points seem worth highlighting.
\begin{alphnum}
\item  The conclusion does not imply that the bounded approximate diagonals for each $A_{i(n)}$ can be assembled to provide a bounded approximate diagonal for $A$ itself. For although $\bigcup_{n\geq 1} A_{i(n)}$ is a subalgebra of $A$, it need not be norm dense in~$A$. Indeed, if $A$ is not unital, then there could exist a non-zero $b\in A$ such that $a_nb=0=ba_n$ for all $a_n\in A_{i(n)}$.
\item Without further information or restrictions on the directed system $(A_i)_{i\in I}$, it does not seem possible to directly relate the amenability constant of $A_{i(n+1)}$ with that of $A_{i(n)}$. For certain nested inclusions of Fourier algebras, it is possible to say more, as will be seen in the proof of Theorem~\ref{t:AMA-discrete-countable}.
\end{alphnum}
\end{rem}

\begin{proof}[The proof of Theorem \ref{t:AMA-discrete-countable}]
Let $\cS$ be the set of countable subgroups of $G$. Given $f\in \FA(G)$, note that $f$ has countable support, since $\FA(G)\subseteq c_0(G)$ and $G$ is discrete, and hence it belongs to $\iota_H(\FA(H))$ for some $H\in\cS$, by Lemma~\ref{l:embed and restrict}.
Thus $\FA(G)=\bigcup_{H\in\cS} \iota_H(\FA(H))$ (we do not need to take closures).
Moreover, $\cS$ is a directed set when ordered by inclusion, and $(\iota_H(\FA(H)))_{H\in \cS}$ is a directed system of closed subalgebras of~$\FA(G)$.

By Proposition \ref{p:AM upper bound sequential}, there is an increasing sequence $H_1 \leq H_2 \leq \dots$ in $\cS$ such that $\AM(\FA(G)) \leq \sup_n \AM(\FA(H_n))$.
Put $\Gamma=\bigcup_{n\geq 1} H_n \in \cS$. Since the restriction maps $\FA(G)\to\FA(\Gamma)$ and $\FA(\Gamma)\to\FA(H_n)$ are contractive  algebra homomorphisms, we have $\AM(\FA(G))\geq \AM(\FA(\Gamma))\geq \AM(\FA(H_n))$ for all $n$.
The result now follows.
\end{proof}

\begin{proof}[The proof of Corollary \ref{c:conjecture-discrete-countable}]
Suppose that there is a discrete counterexample to Conjecture~\ref{conj:main}, i.e.\ a discrete group $G$ such that $\AM(\FA(G))> \AD(G)$. By Theorem \ref{t:AMA-discrete-countable} there is a countable subgroup $\Gamma \leq G$ satisfying
$\AM(\FA(\Gamma)) = \AM(\FA(G)) > \AD(G)$.
Since $\Gamma\leq G$, we have $\AD(G)\geq\AD(\Gamma)$ by \cite[Proposition 2.9]{Choi-minorant}. Therefore $\Gamma$ is also a counterexample to Conjecture \ref{conj:main}.
\end{proof}

\section{A sharper upper bound for the amenability constant}\label{s:prove upper bound}

\subsection{The statement of our main upper bound}

Our goal is to prove the following result.
\begin{thm}[A sharper upper bound for $\AM(\FA(G))$ in the discrete case]
\label{t:sharper upper bound}
Let $G$ be a discrete group that is virtually abelian.
Then
\[ \AM(\FA(G)) \leq 1+ (\maxdeg(G)-1) \left( 1- \frac{1}{|[G,G]|} \right).
\]
\end{thm}

To see where this upper bound comes from, we need to  consider the (standard) method for constructing bounded approximate diagonals for Fourier algebras of discrete groups.
Given a discrete group $G$, let $\Delta\defeq\{(s,s) \colon s\in G\}$ be the diagonal subgroup of $G\times G$.
If $u\in\FA(\Delta)$, then by Lemma \ref{l:embed and restrict} we may view $u$ as an element of $\FA(G\times G)$ that is supported inside $\Delta$, with the same norm. Moreover, since $\Delta\cong G$ we have an isometric algebra isomorphism $\FA(G)\cong\FA(\Delta)$.
We thus obtain an isometric algebra homomorphism $\FA(G)\to\FA(G\times G)$, which (by slight abuse of notation) we denote by $\iota_\Delta$. Explicitly, if $f\in\FA(G)$ then
\[
\iota_\Delta f(s,t)= \begin{cases} f(s) & \text{if $s=t$} \\ 0 & \text{otherwise.}
\end{cases}
\]

For any locally compact group $G$, there is a canonical contractive linear map
\[
J:\FA(G)\ptp\FA(G) \to \FA(G\times G),
\]
defined on elementary tensors by $J(a\otimes b)(s,t) = a(s)b(t)$. If $G$ is also virtually abelian, then by \cite[Theorem 1]{Losert-tp-AG} $J$ is bijective, so that $J^{-1}:\FA(G\times G)\to \FA(G)\ptp\FA(G)$ is well-defined and continuous.

\begin{rem}
When $G$ is virtually abelian, \cite[Remark (a)]{Losert-tp-AG} provides the explicit upper bound $\norm{J^{-1}}\leq \maxdeg(G)$.
An alternative proof of this estimate is given in \cite[Corollary~2.6]{Runde-2006}, but it relies on the fact that $\FA(G\times G)$ coincides with the \emph{operator-space} projective tensor square of $\FA(G)$.
(Warning: in \cite{Runde-2006}, $\ptp$ denotes the operator-space projective tensor product, not the Banach-space projective tensor product.)
\end{rem}

\begin{lem}\label{l:AM via diagonal embedding}
Let $G$ be discrete and virtually abelian. Then $\AM(\FA(G))\leq \norm{J^{-1}\iota_\Delta}$.
\end{lem}

This result may be folklore (it was certainly known to the authors of \cite{Forrest-Runde-2005}). For the reader's convenience, we provide the details.

\begin{proof}[Proof of Lemma \ref{l:AM via diagonal embedding}]
We start by noting that $\iota_\Delta:\FA(G)\to\FA(G\times G)$ is an $\FA(G)$-bimodule map. In particular, since $\FA(G)$ is commutative, $a\cdot\iota_\Delta(f)=\iota_\Delta(f)\cdot a$ for all $a,f\in \FA(G)$.

Since $G$ is virtually abelian it is amenable, so by Leptin's theorem $\FA(G)$ has a bounded approximate identity of norm~$1$, $(f_\alpha)$~say. Put $m_\alpha = J^{-1}\iota_\Delta(f_\alpha) \in \FA(G)\ptp\FA(G)$.

Since $J$ is a bijective $\FA(G)$-bimodule map, so is $J^{-1}$. Hence,
$a\cdot m_\alpha=m_\alpha\cdot a$ for all $a\in\FA(G)$.
Moreover, since $\Delta J^{-1}\iota_\Delta$ is just the identity map on $\FA(G)$, we have $\Delta(m_\alpha)=f_\alpha$.
Thus $(m_\alpha)$ is a bounded approximate diagonal for $\FA(G)$, with norm bounded by $\norm{J^{-1}\iota_\Delta}$.
\end{proof}

\begin{prop}\label{p:bound on diagonal embedding}
Let $G$ be a discrete virtually abelian group. Then
\[
\norm{ J^{-1}\iota_\Delta: \FA(G) \to \FA(G)\ptp\FA(G)} \leq 1 + (\maxdeg(G)-1) \left( 1- \frac{1}{|[G,G]|}\right) ,
\]
with the convention that if $[G,G]$ is infinite then $| [G,G]|^{-1}=0$.
\end{prop}

Note that Theorem \ref{t:sharper upper bound} follows instantly from combining Lemma \ref{l:AM via diagonal embedding} with Proposition~\ref{p:bound on diagonal embedding}, and so it only remains to prove the proposition.
This will be done in two stages. In the first stage, we perform a reduction step that allows us to restrict attention to \emph{countable} virtually abelian groups.
In the second stage, we use the Plancherel theorem for countable virtually abelian groups, and the accompanying formulas for the $\FA(G\times G)$-norm and $\FA(G)\ptp\FA(G)$-norm of a given function in terms of its matrix-valued Fourier coefficients.

\subsection{Reducing to the countable case}
We introduce some notation, to simplify the calculations in this subsection.
For a discrete group~$G$, define
\[ C_G \defeq 1+ (\maxdeg(G)-1) \left(1 - \frac{1}{|[G,G]|}\right) \in [1,\infty]. \]
This is the claimed upper bound on the norm of $J^{-1}\iota_\Delta : \FA(G)\to\FA(G)\ptp\FA(G)$.

\begin{lem}\label{l:monotone}
Let $G$ be discrete and satisfy $\maxdeg(G)<\infty$. Then for every subgroup $L\leq G$, we have
\begin{alphnum}
\item\label{li:monotone maxdeg}
 $\maxdeg(L)\leq\maxdeg(G)$;
\item\label{li:monotone upper bound}
 $C_L\leq C_G$.
\end{alphnum}
\end{lem}

\begin{proof}
Part \ref{li:monotone maxdeg} appears to be folklore; a proof can be given by following the arguments used in the proof of \cite[Proposition 3.1]{Moore-fd-irrep}.
For the reader's convenience, we have included details in the appendix (Section~\ref{app:maxdeg}).

Combining part \ref{li:monotone maxdeg} with the trivial inequality $|[L,L]| \leq | [G,G]|$ yields
\[
\begin{aligned}
C_L-1
& = (\maxdeg(L)-1) \left( 1- \frac{1}{|[L,L]|}\right) \\
& \leq (\maxdeg(G)-1) \left( 1- \frac{1}{|[G,G]|}\right) 
= C_G -1\,,
\end{aligned}
\]
which proves part \ref{li:monotone upper bound}.
\end{proof}

\begin{lem}\label{l:commuting diagram}
Let $G$ be any discrete group and $H\leq G$ any subgroup; let $\iota_H:\FA(H)\to\FA(G)$ be the canonical embedding. Then the following diagram
commutes.
\[
\begin{tikzcd}
	{\FA(G)} && {\FA(G\times G)} && {\FA(G)\ptp \FA(G)} \\
	{} \\
	{\FA(H)} && {\FA(H\times H)} && {\FA(H)\ptp \FA(H)}
	\arrow["{\iota_\Delta}", from=1-1, to=1-3]
	\arrow["J"', from=1-5, to=1-3]
	\arrow["{\iota_H}", from=3-1, to=1-1]
	\arrow["{\iota_\Delta}"', from=3-1, to=3-3]
	\arrow["{\iota_{H\times H}}"', from=3-3, to=1-3]
	\arrow["{\iota_H \ptp\iota_H}"', from=3-5, to=1-5]
	\arrow["J", from=3-5, to=3-3]
        \end{tikzcd}\]
\end{lem}

\begin{proof}
Since
\[
J(\iota_H\ptp\iota_H)(\delta_s\otimes \delta_t) = \delta_{(s,t)}=\iota_{H\times H}J(\delta_s\otimes\delta_t)
\quad\text{for all $s,t\in H$,}
\]
it follows by linearity and continuity that $J(\iota_H\ptp\iota_H)=\iota_{H\times H}J$. Thus the right-hand square commutes.
A similar argument shows that the left-hand square commutes (alternatively, use the characterization of the ranges of $\iota_{H}$ and $\iota_{H\times H}$ given in Lemma \ref{l:embed and restrict}).
\end{proof}

\begin{proof}[Reduction of Proposition \ref{p:bound on diagonal embedding} to the countable case]
Suppose that for every countable virtually abelian group $H$, we have $\norm{J^{-1}\iota_\Delta: \FA(H)\to \FA(H)\ptp\FA(H)} \leq C_H$.

Now let $G$ be a discrete virtually abelian group; recall that this implies $\maxdeg(G)<\infty$ (Remark \ref{r:LCVA implies Type I}). Moreover, since $\FA(G)\subseteq c_0(G)$ and $G$ is discrete, every element of $\FA(G)$ has countable support.
Choose a sequence $(f_n)$ in $\FA(G)$ such that
\[
\norm{f_n} = 1 \text{ and } \norm{J^{-1}\iota_\Delta(f_n)}\geq \norm{J^{-1}\iota_\Delta :\FA(G)\to\FA(G)\ptp\FA(G)} - 2^{-n} \text{ for all $n$.}
\]
Let $L$ be the subgroup of $G$ generated by $\bigcup_n \supp(f_n)$. Then  $L$ is countable,  and it is virtually abelian (since $G$ is). By Lemma \ref{l:embed and restrict}, $f_n$ can be viewed as a norm-$1$ element of $\FA(L)$. Then, by Lemma \ref{l:commuting diagram}, we have
\[ \norm{ J^{-1}\iota_\Delta(f_n)}_{\FA(G)\ptp\FA(G)} \leq
\norm{ J^{-1}\iota_\Delta(f_n)}_{\FA(L)\ptp\FA(L)} \leq C_L \,, \]
with the final inequality holding because of our starting assumption. Taking the supremum over all $n$ and using Lemma \ref{l:monotone}\ref{li:monotone upper bound} gives
\[
\norm{ J^{-1}\iota_\Delta :\FA(G) \to  \FA(G)\ptp\FA(G)} \leq C_L\leq C_G\,,
\]
as required.
\end{proof}

This ``reduction to the countable case'' is done mostly for technical convenience: it allows us to invoke results from non-abelian Fourier analysis that are usually stated for 2nd-countable unimodular Type~I groups, without needing to discuss versions of the Plancherel theorem for uncountable discrete virtually abelian groups. However, we believe that the techniques used in this reduction step have some independent interest.
In any case, the main examples to which we will apply Theorem \ref{t:sharper upper bound} are countable.

\subsection{The proof of Proposition \ref{p:bound on diagonal embedding} for countable groups}
Recall that discrete virtually abelian groups are Type~I (see Remark \ref{r:LCVA implies Type I}).

Throughout this subsection, we adopt the following conventions \textbf{unless explicitly stated otherwise}: $G$~denotes a countable virtually abelian group; $\nu$ denotes Plancherel measure on~$\widehat{G}$, normalized so that $1= \int_{\widehat{G}} d_\pi \,d\nu(\pi)$; $\Omega$ denotes the set of $1$-dimensional representations of~$G$; and we write $d$ for $\maxdeg(G)$.

As in \cite{Runde-2006}, a key role is played by the following estimate, which is purely finite-dimensional in nature.
For a proof, see {\cite[Lemma]{Losert-tp-AG}} or \cite[Lemma~3.6.1]{Kaniuth-Lau-book}.

\begin{lem}\label{l:tp trace-class}
Given $n\in\Nat$, let $S^1(\Cplx^n)$ denote the space of $n\times n$ matrices over $\Cplx$, equipped with the Schatten-$1$ norm. Then for any $m,n\in\Nat$, the natural map $J_{m,n}:S^1(\Cplx^m)\ptp S^1(\Cplx^n)\to S^1(\Cplx^{mn})$ satisfies $\norm{(J_{m,n})^{-1}} = \min(m,n)$.
\end{lem}

(Strictly speaking, the statement of \cite[Lemma]{Losert-tp-AG} concerns the adjoint map $\Phi_{m,n}=(J_{m,n})^\ast$, but it is a standard fact that the norm of a linear map between Banach spaces is equal to the norm of the adjoint map between the duals of these spaces.)

To simplify the formulas in the following calculations, define a function $\Psi: c_{00}(G\times G)\to [0,\infty)$  by
\begin{equation}\label{eq:auxiliary}
\Psi(u) \defeq \int_{\chi\in\Omega} \int_{\sigma\in\widehat{G}} \norm{ (\chi\otimes\sigma)(u)}_{S^1(\sigma)} \,d\nu(\sigma)\,d\nu(\chi)\,.
\end{equation}

\begin{rem}\label{r:absorb}
The notation in Equation \eqref{eq:auxiliary} may require some explanation. If $H_\chi$ and $H_\sigma$ are the Hilbert spaces on which we represent $\chi$ and $\sigma$, then  $(\chi\otimes\sigma)(u)$ is an operator on the Hilbert space $H_\chi\otimes^2 H_\sigma$. But since $H_\chi$ is one-dimensional, $H_\chi\otimes^2 H_\sigma$ may be identified isometrically with~$H_\sigma$. Moreover, using the same convention for Schatten-classes as in Equation \eqref{eq:A(G)-norm via plancherel}, we have \emph{isometric} identifications
$S^1(\chi\otimes\sigma) = S^1(\sigma) = S^1(\chi)\ptp S^1(\sigma)$.
\end{rem}

\begin{lem}\label{l:small saving}
Let $u\in c_{00}(G\times G)\subseteq \FA(G\times G)$. Then
\[
\norm{J^{-1}u}_{\FA(G)\ptp\FA(G)} \leq d\norm{u}_{\FA(G\times G)} - (d-1) \Psi(u).
\]
\end{lem}

\begin{proof}
We use the formulas in Proposition \ref{p:norms via FT} (see Remark \ref{r:scope} for why these formulas are valid for elements of $c_{00}(G\times G)$).

Using Fubini's theorem to write $\int_{\widehat{G}\times\widehat{G}}$ as an iterated integral, and then splitting $\int_{\widehat{G}} \int_{\widehat{G}} = \int_{\Omega}\int_{\widehat{G}} + \int_{\widehat{G}\setminus\Omega}\int_{\widehat{G}}$, we obtain
\[
\norm{u}_{\FA(G\times G)}  = \Psi(u) + \int_{\pi\in\widehat{G}\setminus\Omega} \int_{\sigma\in\widehat{G}} \norm{ (\pi\otimes\sigma)(u)}_{S^1(\pi \otimes\sigma)} \,d\nu(\sigma)\,d\nu(\pi)\;,
\]
and
\[
\norm{J^{-1}u}_{\FA(G)\ptp\FA(G)}  = \Psi(u) + \int_{\pi\in\widehat{G}\setminus\Omega} \int_{\sigma\in\widehat{G}} \norm{ (\pi\otimes\sigma)(u)}_{S^1(\pi)\ptp S^1(\sigma)} \,d\nu(\sigma)\,d\nu(\pi)\;.
\]
Applying Lemma \ref{l:tp trace-class} to the second term in each equation, we deduce that
\[
\norm{J^{-1}u}_{\FA(G)\ptp\FA(G)}  - \Psi(u)
\leq d\left( \norm{u}_{\FA(G\times G)}  - \Psi(u)\right),
\]
and rearranging gives the desired result.
\end{proof}

\begin{lem}\label{l:key trick}
Let $f\in c_{00}(G)$. Then $\Psi(\iota_\Delta(f)) = \nu(\Omega)\ \norm{f}_{\FA(G)}$.
\end{lem}

\begin{proof}
By definition,
\[
\Psi(\iota_\Delta(f)) = \int_{\chi\in\Omega} \int_{\sigma\in\widehat{G}} \norm{ (\chi\otimes\sigma)(\iota_\Delta(f))}_{S^1(\sigma)} \,d\nu(\sigma)\,d\nu(\chi)\;.
\]
But for each $\sigma\in \widehat{G}$ and each $\chi\in\Omega$, we have
\[
(\chi\otimes\sigma)(\iota_\Delta(f)) = \sum_{s,t\in G} \iota_\Delta(f)(s,t) \chi(s)\sigma(t) = \sum_{t\in G} f(t)\chi(t) \sigma(t) = \sigma(f\chi)\;,
\]
and so
\[ 
\Psi(\iota_\Delta(f))
 =\int_{\chi\in\Omega} \int_{\sigma\in\widehat{G}} \norm{ \sigma(f\chi)}_{S^1(\sigma)} \,d\nu(\sigma)\,d\nu(\chi) 
\\
 =\int_{\chi\in\Omega} \norm{  f\chi }_{\FA(G)}\;d\nu(\chi)\;,
\]
where we used Proposition \ref{p:norms via FT}\ref{li:AG-norm} on the inner integral.

We now recall a general result about Fourier algebras (valid for all locally compact groups). Let $\FS(G)$ denote the Fourier--Stieltjes algebra of $G$: then for each $f\in\FA(G)$ and $h\in \FS(G)$ we have $\norm{fh}_{\FA(G)}\leq \norm{f}_{\FA(G)}\norm{h}_{\FS(G)}$. Applying this with $h=\chi$ and $h=\overline{\chi}$, for some given $\chi\in\Omega$, yields
\[
\norm{f\chi}_{\FA(G)}\leq \norm{f}_{\FA(G)} = \norm{f\chi\overline{\chi}}_{\FA(G)} \leq \norm{f\chi}_{\FA(G)}\;.
\]
Thus, for each $\chi\in\Omega$, multiplication by $\chi$ acts as an \emph{isometry} on $\FA(G)$. Hence
\[
\Psi(\iota_\Delta(f))
=\int_{\chi\in\Omega} \norm{  f }_{\FA(G)}\;d\nu(\chi) = \nu(\Omega)\ \norm{f}_{\FA(G)}\,,
\]
as required.
\end{proof}

Our final ingredient, which explains the presence of $|[G,G]|$ in Proposition \ref{p:bound on diagonal embedding}, has nothing to do with amenability constants or Fourier algebras. We suspect that it is not a new result.

\begin{prop}[The Plancherel measure of the set of 1-dimensional representations]\label{p:nu(Omega) and [G,G]}
Let $G$ be a discrete group.
\begin{alphnum}
\item\label{li:[G,G] infinite}
If $[G,G]$ is infinite, then $\nu(\Omega)=0$.
\item\label{li:[G,G] finite}
If $[G,G]$ is finite, then $\nu(\Omega) = |[G,G]|^{-1}$.
\end{alphnum}
\end{prop}


\begin{proof}
Part \ref{li:[G,G] infinite} was proved in \cite[Lemma 5.1]{Choi-minorant}, so we only need to prove part \ref{li:[G,G] finite}.
Suppose that $[G,G]$ is finite, and define $h\in c_{00}(G)$ by
\[
h(t) = \begin{cases}
|[G,G]|^{-1} & \text{if $t\in [G,G]$,} \\
0 & \text{otherwise.}
\end{cases}
\]
Taking $K=[G,G]$ in Proposition \ref{p:pi of mu_K}, the measure $\mu_K$ corresponds to the function $h$ defined above, while the set of $\pi\in\widehat{G}$ that factor through $G \to G / [G,G]$ is exactly the set $\Omega$. We thus obtain $\chi(h) = 1$ for all $\chi\in\Omega$ and $\pi(h)=0$ for all $\pi\in\widehat{G}\setminus\Omega$.
Applying Equation \eqref{eq:plancherel formula} (the Plancherel formula) now gives
\[
\frac{1}{|[G,G]|} = \norm{h}_2^2
= \int_{\widehat{G}} \Tr(\pi(h)\pi(h)^\ast) \,d\nu(\pi) = \int_{\Omega} |\chi(h)|^2\,d\nu(\chi) = \nu(\Omega),
\]
as required.
\end{proof}

\begin{rem}
It was already shown in \cite[Section~4]{Choi-minorant} that $|[G,G]|=2 \iff \nu(\Omega)=1/2$, and this played a crucial role in characterizing those non-abelian groups $G$ on which $\AD$ is minimized.
However, the methods used in that paper do not seem to generalize to larger values of $|[G,G]|$.
\end{rem}

\begin{proof}[The proof of Proposition \ref{p:bound on diagonal embedding} for countable virtually abelian groups]
Let $f\in c_{00}(G)$.
Combining Lemmas \ref{l:small saving} and \ref{l:key trick}, and recalling that $\norm{\iota_\Delta(f)}_{\FA(G\times G)} = \norm{f}_{\FA(G)}$,
 we obtain
\begin{equation}\label{eq:stepping stone}
\norm{J^{-1}\iota_\Delta(f)} \leq d\norm{f}_{\FA(G)} - (d-1)\nu(\Omega)\norm{f}_{\FA(G)}\;.
\tag{$\ast$}
\end{equation}
Since $c_{00}(G)$ is dense in $\FA(G)$, we deduce that
\[ \norm{J^{-1}\iota_\Delta} \leq d -(d-1)\nu(\Omega) = 1+ (d-1)(1-\nu(\Omega)) . \]
Applying Proposition~\ref{p:nu(Omega) and [G,G]} completes the proof.
\end{proof}

\subsection{Applications to discrete groups with two degrees of irreps}\label{s:applications}

\begin{proof}[Proof of Theorem \ref{t:two degrees}]
Let $G$ be a countable non-abelian group, let $d\in\Nat$, and suppose that $\{ d_\pi \colon \pi\in\widehat{G}\} = \{1,d\}$. (This implies, in particular, that $G$ is virtually abelian.) Let $\nu$ be the Plancherel measure for $G$, normalized so that
$1 = \int_{\widehat{G}} d_\pi\,d\nu(\pi)$.

From Theorem~\ref{t:sharper upper bound} we have
\begin{equation}\label{eq:two degrees upper bound}
\AM(\FA(G)) \leq 1 + (d-1)\left( 1- \frac{1}{|[G,G]|}\right).
\tag{$*$}
\end{equation}
We now seek a matching lower bound on $\AD(G)$.
Let $\Omega_n=\{ \pi\in\widehat{G}\colon d_\pi=n\}$. By assumption, $\widehat{G}=\Omega_1\cup \Omega_d$, and so by \cite[Theorem 1.5]{Choi-minorant}
\begin{equation}\label{eq:AD for two degrees}
\AD(G) = \int_{\widehat{G}} (d_\pi)^2\,d\nu(\pi) 
 = \nu(\Omega_1) + \nu(\Omega_d)d^2 \,.
 \tag{$**$}
\end{equation}
However, by our normalization of $\nu$, we have $1 =\nu(\Omega_1) + \nu(\Omega_d)d$. Combining this with \eqref{eq:AD for two degrees} gives
\[ 
\AD(G) = \nu(\Omega_1) + ( 1- \nu(\Omega_1))d = 1 + (d-1)(1-\nu(\Omega_1)),
\]
and using Proposition \ref{p:nu(Omega) and [G,G]}, we deduce that
\begin{equation}\label{eq:two degrees lower bound}
\AD(G) = 1 + (d-1)\left( 1- \frac{1}{|[G,G]|}\right).
\tag{$***$}
\end{equation}
Since $\AM(\FA(G))\geq \AD(G)$, we deduce from \eqref{eq:two degrees upper bound} and \eqref{eq:two degrees lower bound} that
\[
\AM(\FA(G)) = \AD(G) = 1 + (d-1)\left( 1- \frac{1}{|[G,G]|}\right),
\]
as required.
\end{proof}

Our second application is to discrete groups $G$ for which $|G\mathbin{:}Z(G)|=4$. For context, we recall the results stated in Theorem~\ref{t:minimal AM(A(non-abelian))}\ref{li:AD geq 3/2}
and
Theorem~\ref{t:minimal AM(A(non-abelian))}\ref{li:minimal AD}, which were already proved in \cite{Choi-minorant}. Namely, given a locally compact non-abelian group $G$, we always have $\AM(\FA(G))\geq \AD(G)\geq 3/2$, and the second inequality is an equality precisely when $|G\mathbin{:}Z(G)|=4$.

\begin{proof}[The proof of Theorem \ref{t:minimal AM(A(non-abelian))}\ref{li:partial converse}]
Let $G$ be discrete and satisfy $|G\mathbin{:}Z(G)|=4$. In particular, $G$ is non-abelian, so as remarked above, \cite[Theorem 1.6]{Choi-minorant} implies that $\AM(\FA(G))\geq 3/2$. We need to prove the converse inequality.

It is folklore that if $|G\mathbin{:}Z(G)|=4$ then $|[G,G]|=2$ and that each $\pi\in\widehat{G}$ satisfies $d_\pi\in\{1,2\}$.(For details, see the proof of \cite[Theorem 4.6]{Choi-minorant} and combine it with \cite[Lemma 4.2(i)]{Choi-minorant}.)
Substituting these values into Theorem \ref{t:sharper upper bound} gives $\AM(\FA(G)) \leq 3/2$, as required.
\end{proof}

\section{New examples of groups satisfying the main conjecture}\label{s:p-reduced Heis}
Given a ring $S$ (assumed to be unital and commutative), we define
\begin{equation}
\sH(S) \defeq \left\{ \utmat{x}{z}{y} \colon x,y,z\in S\right\}
\end{equation}
which forms a group with respect to matrix multiplication. We refer to $\sH(S)$ as the \dt{Heisenberg group over $S$}, or the $S$-Heisenberg group. It is convenient to view elements of $\sH(S)$ as ordered triples $(x,y,z)$, with the group multiplication given by the formula
\begin{equation}\label{eq:heisenberg law}
(x,y,z)(x',y',z') \defeq (x+x', y+y', z+xy'+z')
\end{equation}
and we will adopt this notation going forwards.

\begin{rem}\label{r:Heisenberg as semidirect}
For studying the representation theory of $\sH(S)$, it is useful to observe that it decomposes as an internal semi\-direct product, with the subgroup $\{(x,0,0)\colon x\in S\}$ acting on the normal abelian subgroup $\{(0,y,z)\colon y,z\in S\}$ by conjugation:
\[
(x,0,0)(0,y,z)(x,0,0)^{-1} = (x,0,0)(0,y,z)(-x,0,0) = (x,0,0)(-x,y,z) = (0,y,xy+z).
\]
We denote this decomposition by $\sH(S)\cong S\ltimes (S\oplus S)$.
\end{rem}

The following facts are easily verified:
\begin{itemize}
\item the function $(x,y,z)\mapsto (x,y)$ is a group homomorphism from $\sH(S)$ onto the additive group of $S\oplus S$; 
\item the subgroup $\{(0,0,z)\colon z\in S\}$ is contained in the centre of $\sH(S)$, and is isomorphic to $(S,+)$.
\end{itemize}

\begin{dfn}\label{d:p-reduced Heisenberg}
Let $p$ be prime, and assume that $p$ (when viewed as $\overbrace{1_S+\dots + 1_S}^{p} \in S)$ is not invertible in~$S$. We define $\Hpred(S)$ to be the quotient of $\sH(S)$ by the (central) subgroup $\{ (0,0,ps) \colon s\in S\}$.
\end{dfn}

\begin{rem}\label{r:reduced Heisenberg}
The group $\sH(\Real)$ has centre isomorphic to $(\Real,+)$, and quotienting by a central copy of $\bbZ$ gives a Lie group, the so-called (real) \dt{reduced Heisenberg group}, which turns up naturally in connection with Gabor time-frequency analysis. The groups $\Hpred(S)$ may be thought of as a discretized analogue of the ``true'' reduced Heisenberg group.
\end{rem}

The condition $p\notin S^\times$, which is part of our definition of $\Hpred(S)$, is satisfied in both of the following cases: $S=\bbZ$, the usual ring of integers; or $S=\Zadic{p}$, the ring of $p$-adic integers.
It is convenient to establish some general properties of $\Hpred(S)$, allowing us to treat both of these cases simultaneously.

We may identify $\Hpred(S)$ as a set with $S\times S\times (S/pS)$ and denote its elements by $(x,y,[z])$, with $[z]$ denoting the image of $z$ in $S/pS$. With this convention, we have
\[
(1,0,[0])(0,1,[0]) = (1,1,[1]) \text{ and } (0,1,[0])(1,0,[0])=(1,1,[0]),
\]
showing that $\Hpred(S)$ is non-abelian. Note also that the semi\-direct product decomposition $\sH(S)\cong S\ltimes (S\oplus S)$ from Remark \ref{r:Heisenberg as semidirect} induces a semi\-direct product decomposition $\Hpred(S)\cong S\ltimes (S\oplus S/pS)$.

\begin{prop}\label{p:not finite times anything}
Assume that the additive group of $S$ is torsion-free, and suppose $p$ is a prime that is not invertible in~$S$. If there exist groups $B$ and $C$ with $B$ finite, such that $\Hpred(S)\cong B\times C$, then $B$ is trivial.
\end{prop}

\begin{proof}
For this proof, we put $G=\Hpred(S)$ and $K=\{(0,0,[s]) \colon s\in S\}\subset G$, to reduce notational clutter. Observe that if $F$ is any finite subgroup of $G$, then $F\subseteq K$ (because $G/K \cong S\oplus S$ is torsion-free); since $K$ is cyclic of order~$p$, it follows that $F$ is either trivial or all of~$K$.

Now suppose that $G\cong B\times C$ for some groups $B$ and $C$ where $B$ is finite and non-trivial; we will derive a contradiction. The previous paragraph implies that the given isomorphism $B\times C \to G$ maps $B\times\{1_C\}$ bijectively onto~$K$. Hence, $C\cong G/K \cong S\oplus S$, which is abelian. Since $B\cong K$ is also abelian, we conclude that $G$ is abelian, which is not the case.
\end{proof}

\begin{proof}[The proof of Theorem \ref{t:new examples}\ref{li:p-reduced integer Heis}]
The unitary representations of $\Hpred(\bbZ)$ can be worked out using the Mackey machine (see \cite[Section 6.8.2]{Folland-book} for a detailed explanation). In particular, if $\pi\in \widehat{\Hpred(\bbZ)}$ then $d_\pi \in\{1,p\}$.
Moreover, a straightforward calculation shows that the commutator subgroup of $\Hpred(\bbZ)$ has order~$p$. Applying Theorem \ref{t:two degrees}, we deduce that
\[ \AM(\FA(\Hpred(\bbZ)))= \AD(\Hpred(\bbZ)) = p - 1 + \frac{1}{p}\,, \]
as claimed.
\end{proof}

We now turn to the proof of the corresponding statement for $\Hpred(\Zadic{p})$. This group is not discrete, so we cannot apply Theorem \ref{t:two degrees}; instead, our leverage comes from the fact that this group is profinite.
We require some basic facts concerning profinite groups and their natural compact topologies: see
the second half of \cite[Section 1.8]{Deitmar-Echterhoff-book} for a quick sketch, or \cite[Section 1.1]{Ribes-Zalesskii-book} for a more comprehensive account.

The following proposition summarizes some general properties of Fourier algebras of profinite groups, which all seem to be folklore.

\begin{prop}\label{p:A(profinite)}
Let $(I,\preceq)$ be a directed set and let $(G_i)_{i\in I}$ be an inverse system of finite groups. Let $G=\varprojlim_{i\in I} G_i$ equipped with its natural compact topology. Write $Q_i:G\to G_i$ for the canonical projection, and let $Q_i^*:C(G_i) \to C(G)$ be the algebra homomorphism induced by composition with $Q_i$.
Then
\begin{alphnum}
\item\label{li:subalgebra}
 for each $i\in I$, $Q_i^*$ maps $\FA(G_i)$ isometrically onto a norm-closed subalgebra of $\FA(G)$;
\item\label{li:exhaustion}
 $\bigcup_{i\in I} Q_i^*(\FA(G_i))$ is norm-dense in $\FA(G)$;
\item\label{li:upper bound}
$\AM(\FA(G)) \leq \sup_{i\in I} \AM(\FA(G_i))$.
\end{alphnum}
\end{prop}


\begin{proof}
Part \ref{li:subalgebra} is a special case of general results concerning quotient homomorphisms between locally compact groups and the maps that these induce between the respective Fourier(--Stieltjes) algebras. See e.g. \cite[Corollary 2.2.4]{Kaniuth-Lau-book} or \cite[Proposition~2.4.2]{Kaniuth-Lau-book}.

Since we could not find an explicit and self-contained reference for part \ref{li:exhaustion}, we provide a proof here for the reader's convenience.
Consider the following space:
\[
\cF_{00} \defeq \bigcup_i Q_i^\ast (C(G_i)) \subseteq C(G).
\]
It is easily checked that $\cF_{00}$ is a subalgebra of $C(G)$ that is closed under complex conjugation. Moreover, the standard construction of $G$ realizes it as a subset of $\prod_i G_i$
(see \cite[Pro\-pos\-ition 1.8.9]{Deitmar-Echterhoff-book} or \cite[Pro\-pos\-ition 1.1.1]{Ribes-Zalesskii-book});
in particular, given $x,x'\in G$ with $x\neq x'$, there exists some $j\in I$ such that $Q_j(x)\neq Q_j(x')$.
Thus, $\cF_{00}$ separates the points of $G$. By the Stone--Weierstrass theorem, it follows that $\cF_{00}$ is dense in $C(G)$ with respect to the supremum norm, and therefore it is dense in $L^2(G)$ with respect to the $L^2$-norm.
Now let $f\in\FA(G)$. Since $f=\alpha_G(\xi\otimes\eta)$ for some $\xi,\eta\in L^2(G)$, a routine approximation argument produces $\xi_n,\eta_n\in \cF_{00}$ such that $\norm{\alpha_G(\xi_n\otimes\eta_n)-f}_{\FA(G)}\to 0$. By the definition of $\cF_{00}$, each function $\alpha_G(\xi_n\otimes\eta_n)$ belongs to $\bigcup_{i\in I} Q_i^*(\FA(G_i))$, and thus we have proved part~\ref{li:exhaustion}.

Finally, part \ref{li:upper bound} follows immediately from part \ref{li:exhaustion} and Lemma \ref{l:AM of inductive limit}.
\end{proof}

\begin{rem}
One can use Proposition \ref{p:AM upper bound sequential} to give a slight refinement of Proposition \ref{p:A(profinite)}\ref{li:upper bound}: namely, there is an increasing sequence $i(1)\preceq i(2)\preceq \dots$ in $I$, such that
\[ \AM(\FA(G)) \leq \limsup_{n\geq 1} \AM(\FA(G_{i(n)})) \;. \]
\end{rem}

\begin{proof}[The proof of Theorem \ref{t:new examples}\ref{li:p-reduced p-adic Heis}]
To ease notation, let $\bbS_n$ denote the ring $\bbZ/p^n\bbZ$.

By construction, $\Zadic{p}$ can be identified as a topological group with the inverse limit $\varprojlim_n \bbS_n$; this identification also respects the ring structure. It follows, either by direct calculations or the universal property of inverse limits, that $\Hpred(\Zadic{p})$ can be identified as a topological group with the inverse limit $\varprojlim_n \Hpred(\bbS_n)$.

As previously mentioned (see the comments before Remark \ref{r:reduced Heisenberg}), we have a semi\-direct product decomposition
\[
\Hpred(\bbS_n) \cong \bbS_n\ltimes  \left( \bbS_n \oplus  (\bbS_n/p\bbS_n) \right)
\cong \bbS_n\ltimes  \left( \bbS_n \oplus \bbS_1 \right),
\]
and so the irreducible representations of $\Hpred(\bbS_n)$ can be determined (up to equivalence) using Mackey theory. One finds that $\Hpred(\bbS_n)$ has exactly $p^{2n}$ $1$-dimensional representations and exactly $p^{2n-1}-p^{2n-2}$ irreducible representations of degree $p$, and no others.

Thus, for each $n$,
applying Johnson's formula \eqref{eq:johnson formula} yields
\[
\AM(\FA(\Hpred(\bbS_n)))= \frac{1}{p^{2n+1}} \left( p^{2n}\cdot 1^3 + (p^{2n-1}-p^{2n-2})\cdot p^3\right ) = p-1+p^{-1} .
\]
We deduce, using Proposition \ref{p:A(profinite)}\ref{li:upper bound}, that $\AM(\FA(\Hpred(\Zadic{p})))\leq p-1+p^{-1}$.

Since $\AM(\FA(\cdot))\geq \AD(\cdot)$, to complete the proof of the theorem it suffices to prove that $\AD(\Hpred(\Zadic{p}))\geq p-1+p^{-1}$.
The canonical injective ring homomorphism $\bbZ\to\bbZ_p$ induces an injective group homomorphism $\sH(\bbZ)\to\sH(\Zadic{p})$, and a little book-keeping shows that this descends to an injective group homomorphism $\Hpred(\bbZ)\to \Hpred(\Zadic{p})$. We already showed that $\AD(\Hpred(\bbZ)) = p-1+p^{-1}$, and so by monotonicity of $\AD$ with respect to subgroup inclusion (Proposition~\ref{p:AD properties}\ref{li:AD of subgroup}), we conclude that
\[
\AD(\Hpred(\Zadic{p})) \geq \AD(\Hpred(\bbZ))=p-1+p^{-1} \;,
\]
as required.
\end{proof}

Finally, we show that neither $\Hpred(\Zahl)$ nor $\Hpred(\Zadic{p})$ can be written as the product of a finite group with an AD-maximal group.

\begin{proof}[The proof of Theorem \ref{t:new examples}\ref{li:not cheap} ]
Let $G$ be either $\Hpred(\Zahl)$ or $\Hpred(\Zadic{p})$, and suppose that $G\cong F\times H$ where $F$ is finite and $H$ is AD-maximal. By Proposition \ref{p:not finite times anything} $F$ has to be trivial, so $G=H$ is AD-maximal. This is a contradiction, since $\maxdeg$ takes positive integer values while $\AD(G)=p-1+p^{-1}$ is not an integer.
\end{proof}

\section{Concluding remarks and questions}\label{s:where next}
The results in this paper provide further evidence for Conjecture \ref{conj:main}, and thus provide further evidence for a positive answer to Question~\ref{q:is AMAG multiplicative}. Nevertheless, that question remains open. We draw attention to a special case that might be more accessible.

\begin{qn}
Let $G$ be a discrete group and let $F$ be a finite group. Is $\AM(\FA(G\times F))$ equal to $\AM(\FA(G))\ \AM(\FA(F))$?
\end{qn}

Since $\AD(G)$ only depends on the underlying group of $G$, Conjecture \ref{conj:main} would also imply a positive answer to the following question.

\begin{qn}\label{q:ignore topology}
Let $G$ be a locally compact group and let $G_d$ denote the same group equipped with the discrete topology. Is $\AM(\FA(G))$ equal to $\AM(\FA(G_d))$?
\end{qn}

Recall that Theorem \ref{t:sharper upper bound} provides a sharper upper bound on $\AM(\FA(G))$ for $G$ discrete and virtually abelian. Our proof of the theorem relies on the fact that the Plancherel measure for such $G$ is finite, but Plancherel measure does not appear in the conclusion of the theorem. It is therefore natural to ask the following question.

\begin{qn}\label{q:general upper bound}
Can Theorem \ref{t:sharper upper bound} be extended to cover all locally compact virtually abelian groups? (The question is only interesting if $[G,G]$ is finite.)
\end{qn}

A positive answer to Question \ref{q:general upper bound} would lead to a stronger version of Theorem \ref{t:minimal AM(A(non-abelian))}\ref{li:partial converse}, in which we do not need to assume that $G$ is discrete. This would answer \cite[Question~6.2]{Choi-minorant} in full generality.

Finally, recall that in order to confirm Conjecture \ref{conj:main} for the groups $\Hpred(\Zadic{p})$, we exploited their structure as profinite groups (not just as compact groups). A natural avenue for further work would be to find more general classes of profinite groups $G$ for which $\AM(\FA(G))$ can be explicitly calculated.
Progress in this direction would be aided by a positive answer to the following question.

\begin{qn}
Is the upper bound in Proposition \ref{p:A(profinite)}\ref{li:upper bound} an equality, possibly after passing to a subsystem?
\end{qn}

\section*{Acknowledgments}
The work presented here began through discussions between the authors during a visit to Winnipeg in 2023, following a meeting of the \textit{Canadian Abstract Harmonic Analysis Symposium},
and was further developed in March 2024 during the Research In Teams workshop \textit{Homological invariants of Fourier algebras} (24rit022) at the Banff International Research Station (BIRS).
Important progress towards Theorem \ref{t:sharper upper bound} was made through discussions during the \textit{Banach Algebras and Operator Algebras 2024} conference at the University of Waterloo, Canada. 
Many technical details, including the proof of Theorem \ref{t:AMA-discrete-countable} and most of Section \ref{s:p-reduced Heis}, were worked out in April 2025 during the Research in Groups workshop \textit{Amenability constants of Fourier algebras}, held at the International Centre for Mathematical Sciences (ICMS) in Edinburgh.
The authors thank the respective organizers and hosts of these workshops for providing stimulating working environments, and generous help with local arrangements.
MG~acknowledges financial support through Simons Travel Support for Mathematicians
and National Science Foundation grant \href{https://www.nsf.gov/awardsearch/show-award/?AWD_ID=2408008&HistoricalAwards=false}{DMS-2408008}.
YC~acknowledges practical support from N.~J.~Laustsen and D.~Pauksztello during a period of severely reduced mobility.

\appendix

\section{Additional proofs}

\subsection{Properties of the $\maxdeg$ function}\label{app:maxdeg}

In Remark \ref{r:products of AD-maximal}, we mentioned that the product of two AD-maximal groups $G_1$ and $G_2$ is itself AD-maximal. We now give a brief explanation. First note that $G_1\times G_2$ is virtually abelian (since both $G_1$ and $G_2$ are). Then, since AD is multiplicative (Proposition \ref{p:AD properties}\ref{li:AD of product}),
it suffices to show that $\maxdeg(G_1\times G_2)=\maxdeg(G_1)\maxdeg(G_2)$.
This will follow once we know that the natural map from $\widehat{G_1}\times\widehat{G_2}$ to $\widehat{G_1\times G_2}$, given by tensoring representations, is bijective. The quickest way to justify this (although it is overkill) is to recall that since $G_1$ is virtually abelian it is Type~I (see Remark \ref{r:LCVA implies Type I}), and then appeal to \cite[Theorem 7.17]{Folland-book}.

We also stated, in Lemma~\ref{l:monotone}\ref{li:monotone maxdeg}, that for any discrete group $G$ with $\maxdeg(G)<\infty$ and any subgroup $L\leq G$ we have $\maxdeg(L)\leq\maxdeg(G)$. For the reader's convenience, we provide a proof of this fact. Our argument is adapted from the proof of \cite[Proposition~3.1]{Moore-fd-irrep}, and in particular we use the following result of Amitsur and Levitzki~\cite{Amitsur-Levitzki}.

\begin{lem}\label{l:AL identity}
For each $n\in\Nat$, there exists a noncommutative polynomial $P_{2n}$ in $2n$ variables, with the following properties:
\begin{itemize}
\item if $1\leq d\leq n$, then $P_{2n}(X_1,\dots, X_{2n})=0$ for all $X_j\in M_d(\Cplx)$;
\item if $d\geq n+1$, then there exist $Y_1,\dots, Y_{2n}\in M_d(\Cplx)$ such that $P_{2n}(Y_1,\dots, Y_{2n})\neq 0$.
\end{itemize}
We say that a unital ring $A$ \dt{satisfies $P_{2n}$} if $P_{2n}(a_1,\dots, a_{2n})=0$ for all $a_1,\dots, a_{2n}\in A$.
\end{lem}

Put $m=\maxdeg(G)$, and consider the complex group ring $\Cplx G$.
For each $\pi\in\widehat{G}$, the ring $\pi(\Cplx G)$ is contained in $\Bdd(H_\pi)$ where $d_\pi \leq m$, and so it satisfies~$P_{2m}$.
It is a standard but non-trivial fact that $\widehat{G}$ separates the points of $\Cplx G$.
(In fact, for any locally compact group, its irreducible unitary representations separate the points of its $L^1$-convolution algebra: see e.g. \cite[Corollary~7.2]{Folland-book}.)
We deduce that $\Cplx G$ also satisfies~$P_{2m}$. Hence $\Cplx L$, being a unital subring of $\Cplx G$, also satisfies~$P_{2m}$.

However, for each $\sigma\in\widehat{L}$, it follows from Schur's lemma and the bicommutant theorem that $\sigma(\Cplx L)$ is SOT-dense in $\Bdd(H_\sigma)$. Since multiplication in $\Bdd(H_\sigma)$ is separately SOT-continuous, by taking iterated limits it follows that $\Bdd(H_\sigma)$ also satisfies the identity $P_{2m}=0$. This forces $d_\sigma \leq m$, as required.

\subsection{The norm on the Fourier algebra in terms of the operator-valued Fourier transform}
\label{app:norms via FT}

In this section, we explain how the norm formulas in Proposition \ref{p:norms via FT} can be deduced from known results concerning the operator-valued Fourier transform for second-countable unimodular Type~I groups.
The formulas express the fact that maps between certain Banach spaces are isometries, and so we need to introduce the relevant spaces.

Given such a group $G$, fix a choice of Haar measure on $G$, and let $\nu$ be the corresponding Plancherel measure on $\Ghat$. We write $L^2(\Ghat;S^2)$ for the direct integral of Hilbert spaces $\int^{\oplus}_{\Ghat} S^2(\pi)\,d\nu(\pi)$. Elements of $L^2(\Ghat;S^2)$ are (equivalence classes of) measurable fields of Hilbert-Schmidt operators; if $C=(C_\pi)_{\pi\in\Ghat}$ is such an operator field, then
\[ \norm{C}_{L^2(\Ghat;S^2)} \defeq \left( \int_{\Ghat} \norm{C_\pi}_{S^2(\pi)}^2 \;d\nu(\pi) \right)^{1/2} \]
(We are suppressing certain technicalities regarding measurable selections; see \cite[Section 7.4]{Folland-book} for an accessible summary.)

The Plancherel formula \eqref{eq:plancherel formula} is then a consequence of the following more precise result. See e.g. \cite[Section 8.5]{Deitmar-Echterhoff-book} for a concise account and references to a full proof.

\begin{thm}[Plancherel theorem, precise version]
\label{t:plancherel bijection}
There is a \emph{bijective} linear isometry $\cF_2:L^2(G)\to L^2(\Ghat;S^2)$, that satisfies
\[
\cF_2(f)_\pi =\pi(f) \quad\text{for all $f\in (L^1\cap L^2)(G)$ and $\nu$-a.e. $\pi\in\Ghat$.}
\]
\end{thm}

Note that if $f\in L^1(G)$ then $\pi(f)$ is defined unambiguously for each $\pi\in\Ghat$, but if $f\in L^2(G)\setminus L^1(G)$ then
$\cF_2(f)_\pi$ is only defined for $\nu$-a.e.~$\pi\in\Ghat$.

The other space we need to consider is $L^1(\Ghat;S^1)$. This denotes the space of (equivalence classes of) measurable operator fields $(C_\pi)_{\pi\in \Ghat}$ such that $C_\pi \in S^1(\pi)$ for $\nu$-a.e.~$\pi\in\widehat{G}$ and $\int_{\Ghat} \norm{C_\pi}_{S^1(\pi)}\,d\nu(\pi) <\infty$; this is a Banach space with respect to the natural norm.

\begin{thm}[Inverse Fourier transform]\label{t:isometric IFT}
Given $C=(C_\pi)_{\pi\in\Ghat}\in L^1(\Ghat;S^1)$ and $s\in G$, define
\begin{equation}\label{eq:define IFT}
\Psi_G(C)(s) \defeq \int_{\Ghat} \Tr( C_\pi \pi(s^{-1}))\,d\nu(\pi).
\end{equation}
Then $\Psi_G$ is a bijective linear isometry of $L^1(\Ghat;S^1)$ onto $\FA(G)$.
\end{thm}

\begin{proof}
See \cite[Theorem 3.1]{Lipsman-1974} or \cite[Theorem 4.12]{Fuhr-book}.
\end{proof}

Let $\cV_2^G$ denote the linear span of the set $\{ \alpha_G(\xi\otimes\eta) \colon \xi,\eta\in (L^1\cap L^2)(G)\}$. For every $\xi$ and $\eta$ in $(L^1\cap L^2)(G)$, we have
\begin{equation}\label{eq:the key}
\alpha_G(\xi\otimes\eta)(s) = \langle \lambda_s\xi,\overline{\eta}\rangle_{L^2(G)}
= \int_G \xi(s^{-1}t) \eta(t) \,dt = (\eta *\xi^\vee)(s)
\end{equation}
where $\xi^\vee(p)\defeq \xi(p^{-1})$. (Crucially, since $G$ is unimodular, the map $\xi\mapsto\xi^\vee$ is an isometric involution on both $L^1(G)$ and $L^2(G)$.) In particular, $\cV_2^G\subseteq (\FA\cap L^1)(G)$.

\begin{prop}\label{p:identify FC}
Let $f\in \cV_2^G$, and let $C$ denote the operator field $(\pi(f))_{\pi\in\Ghat}$.
Then $C\in L^1(\Ghat;S^1)$, and $\Psi_G(C)=f$.
\end{prop}

\begin{proof}
By linearity, it suffices to verify these statements when $f=\alpha_G(\xi\otimes\eta)$ for some $\xi,\eta\in (L^1\cap L^2)(G)$.
Since $\eta,\xi^\vee\in L^1(G)$, the calculation in Equation \eqref{eq:the key} implies that $\pi(f)=\pi(\eta)\pi(\xi^\vee)$ for all $\pi\in\Ghat$.
But since $\eta,\xi^\vee\in L^2(G)$, it follows from Theorem \ref{t:plancherel bijection} that $C$ is the pointwise product of the operator fields $\cF_2(\eta)$, $\cF_2(\xi^\vee)\in L^2(\Ghat;S^2)$. Hence $C\in L^1(\Ghat;S^1)$.

To finish, we need to show that for each $s\in G$ the following equality holds:
\[
\Psi_G(\cF_2(\eta) \cF_2(\xi^\vee)  )(s)
= \alpha_G(\xi\otimes\eta)(s) = \langle \eta, \lambda_s\overline{\xi}\rangle_{L^2(G)}\,.
\]
This calculation can be found in most accounts of Fourier inversion for~$G$, see e.g.\ the proofs of \cite[Theorem 3.1]{Lipsman-1974} or \cite[Theorem 4.4]{Fuhr-book}. (The trick is to use the definition of $\Psi_G$ to rewrite the left-hand side as the inner product in $L^2(\Ghat;S^2)$ of the operator fields $\cF_2(\eta)$ and $\cF_2(\lambda_s\overline{\xi})$, and then appeal to the ``isometry'' part of Theorem \ref{t:plancherel bijection}.)
\end{proof}

Combining Proposition \ref{p:identify FC} with the isometry part of Theorem \ref{t:isometric IFT} yields
\[
\norm{f}_{\FA(G)} = \int_{\Ghat} \norm{\pi(f)}_{S^1(\pi)} \,d\nu(\pi) \quad\text{for all $f\in \cV_2^G$,}
\]
which is precisely the assertion of Proposition~\ref{p:norms via FT}\ref{li:AG-norm}.

Now let $N$ be another second-countable unimodular Type~I group. Given $w\in \cV_2^G\otimes\cV_2^N$, we have
\[
\norm{w}_{\FA(G)\ptp\FA(N)} = \norm{ (\Psi_G^{-1}\ptp\Psi_N^{-1})(w)}_{L^1(\Ghat;S^1)\ptp L^1(\Nhat;S^1)} \,.
\]
It is a standard result (due to Grothendieck) that, for measure spaces $(\Omega_1,\mu_1)$ and $(\Omega_2,\mu_2)$, the natural map $L^1(\Omega_1,\mu_1)\ptp L^1(\Omega_2,\mu_2) \to L^1(\Omega_1\times\Omega_2;\mu_1\times\mu_2)$ is an isometric isomorphism of Banach spaces. One can build on this to show that the natural map
\[L^1(\Ghat;S^1)\ptp L^1(\Nhat;S^1) \to L^1(\Ghat\times\Nhat;S^1\ptp S^1) \]
is an isometric isomorphism.
The key technical point is that $L^1(\Ghat;S^1)$ can be decomposed as an $\ell^1$-direct sum of spaces of the form $L^1(\Omega_n;S^1(H_n))$ for $n\in\bbN\cup\{\infty\}$, where $H_n$ is a Hilbert space of dimension $n$ for $n\in \Nat$ and $H_\infty\defeq \ell^2$; and each summand can in turn be identified with the projective tensor product $L^1(\Omega_n)\ptp S^1(H_n)$.
We omit the details, which are measure-theoretic rather than anything particular to do with harmonic analysis.

After making this identification: by using Proposition \ref{p:identify FC} for $G$ and for $N$, we obtain
\[
\norm{w}_{\FA(G)\ptp\FA(N)} = \norm{ ((\pi\otimes\sigma)(w) )_{(\pi,\sigma)\in\Ghat\times\Nhat} }_{L^1(\Ghat\times\Nhat;S^1\ptp S^1)}\;.
\]
Writing out the norm on the right-hand side gives the statement in Proposition \ref{p:norms via FT}\ref{li:AG-tp-norm}.

Finally, one would like to say that Proposition \ref{p:norms via FT}\ref{li:AGN-norm} just follows from Proposition \ref{p:norms via FT}\ref{li:AG-norm} by replacing $G$ with $G\times N$. To make this rigorous: if we write $\nu_G$ and $\nu_N$ for the Plancherel measures of $G$ and $N$ respectively, then we need to identify the measure space $(\widehat{G\times N}, \nu_{G\times N})$ with the product of the measure spaces $(\Ghat,\nu_G)$ and $(\Nhat,\nu_N)$.
Although this identification is often taken for granted, we could not find a suitable reference. (See Remark \ref{r:issues} for a discussion of some of the issues.)

We therefore indicate an alternative approach, based on the fact that we can identify the Hilbert-space tensor product $L^2(G)\otimes^2 L^2(N)$ with $L^2(G\times N)$. By similar arguments, one can identify $L^2(\Ghat;S^2)\otimes^2 L^2(\Nhat;S^2)$ with $L^2(\Ghat\times\Nhat;S^2)$, where elements of the second space are measurable operator fields $(T_{(\pi,\sigma)})_{(\pi,\sigma)\in\Ghat\times\Nhat}$ such that $T_{(\pi,\sigma)}\in S^2(\pi)\otimes^2 S^2(\sigma)\equiv S^2(\pi\otimes\sigma)$.
It then follows from Theorem \ref{t:plancherel bijection} that we have a bijective linear isometry
\[ \cG_2^{(G,N)} : L^2(G\times N) \to L^2(\Ghat\times\Nhat;S^2) \]
that satisfies $\cG_2^{(G,N)} (w)_{(\pi,\sigma)} = (\pi\otimes\sigma)(w)$ for all $w\in\cV_2^G\otimes\cV_2^N$.
Moreover, inspecting the proof of Theorem \ref{t:isometric IFT} shows that there is a bijective linear isometry
\[ \Phi_{(G,N)} : L^1(\Ghat\times\Nhat;S^1) \to \FA(G\times N), \]
defined by the obvious two-variable analogue of the formula \eqref{eq:define IFT}. We can now repeat the argument used in the proof of Proposition \ref{p:identify FC}, to show that 
for all $w\in \cV_2^G\otimes\cV_2^N$ we have 
\[
\norm{w}_{\FA(G\times N)} = \norm{ ((\pi\otimes\sigma)(w) )_{(\pi,\sigma)\in\Ghat\times\Nhat} }_{L^1(\Ghat\times\Nhat;S^1)} \;.
\]
Writing out the norm on the right-hand side gives the statement in Proposition \ref{p:norms via FT}\ref{li:AGN-norm}.

\begin{rem}\label{r:issues}
It is known that for $G$ and $N$ in the class being considered, the natural function $\Ghat\times\Nhat\to \widehat{G\times N}$ is a continuous bijection. However, to identify the two sets as measure spaces, we have to ensure that the product $\sigma$-algebra on $\Ghat\times\Nhat$, arising from the Mackey Borel structures on the two factors, coincides with the Mackey Borel structure on $\widehat{G\times N}$. Moreover, even after one has verified this, a separate calculation would be needed to show that $\nu_G\times\nu_N$ satisfies the Plancherel formula for~$G\times N$.
\end{rem}



\vfill

\newcommand{\address}[1]{{\small\sc#1.}}
\newcommand{\email}[1]{\noindent Email: \texttt{#1}}

\noindent
\address{Yemon Choi,
School of Mathematical Sciences,
Lancaster University,
Lancaster LA1 4YF, United Kingdom} 

\email{y.choi1@lancaster.ac.uk}

\noindent
\address{Mahya Ghandehari,
Department of Mathematical Sciences,
University of Delaware,
Newark, Delaware 19716, United States of America}

\email{mahya@udel.edu}

\end{document}